\documentclass[12pt,reqno]{amsart}
\usepackage{amsmath}
\usepackage{amssymb}
\usepackage{amsthm}
\usepackage{mathrsfs}
\usepackage{latexsym}
\usepackage{xcolor}
\usepackage{comment}

\newtheorem{theorem}{Theorem}[section]
\newtheorem{lemma}{Lemma}[section]
\newtheorem{corollary}{Corollary}[section]
\newtheorem{remark}{Remark}[section]
\newtheorem{definition}{Definition}[section]
\newtheorem{proposition}{Proposition}[section]
\newtheorem{example}{Example}[section]
\newtheorem{assumption}{Assumption}[section]
\numberwithin{equation}{section}
\newcommand{\bth}{\begin{theorem}}
\newcommand{\ethe}{\end{theorem}}

\newcommand{\bre}{\begin{remark}}
\newcommand{\ere}{\end{remark}}

\newcommand{\ble}{\begin{lemma}}
\newcommand{\ele}{\end{lemma}}
\newcommand{\bde}{\begin{definition}}
\newcommand{\ede}{\end{definition}}

\newcommand{\bco}{\begin{corollary}}
\newcommand{\eco}{\end{corollary}}

\newcommand{\bpr}{\begin{proposition}}
\newcommand{\epr}{\end{proposition}}

\newcommand{\bexer}{\begin{exercise}}
\newcommand{\eexer}{\end{exercise}}
\newcommand{\breh}{\begin{hint}}
\newcommand{\ereh}{\end{hint}}

\newcommand{\halmos}{\hfill \qed}

\newcommand{\bexam}{\begin{example}}
\newcommand{\eexam}{\end{example}}

\newcommand{\pr} {{\bf Proof.}}

\newcommand{\bfi}{\begin{fig}}
\newcommand{\efi}{\end{fig}}

\newcommand{\beao}{\begin{eqnarray*}}
\newcommand{\eeao}{\end{eqnarray*}\noindent}

\newcommand{\beam}{\begin{eqnarray}}
\newcommand{\eeam}{\end{eqnarray}\noindent}
\newcommand{\E}{\mathbf{E}}
\newcommand{\PP}{\mathbf{P}}

\newcommand{\xto}{x\to\infty}

\newcommand{\bF}{\overline{F}}

\newcommand{\bG}{\overline{G}}
\newcommand{\bV}{\overline{V}}

\newcommand{\bbr}{{\mathbb R}}

\newcommand{\bbb}{{\mathbb B}}

\newcommand{\bbn}{{\mathbb N}}

\newcommand{\vep}{\varepsilon}

\allowdisplaybreaks[1]

\begin{document}
\title[Local uniform asymptotics for a non-standard risk model]{Local uniform asymptotics for a non-standard risk model and interplay of insurance and financial risks stemming by systemic factors}

\author[D.G. Konstantinides, C.D. Passalidis, M. Yuan]{Dimitrios G. Konstantinides, Charalampos  D. Passalidis, Meng Yuan}
\address{Dept. of Statistics and Actuarial-Financial Mathematics,
University of the Aegean, Karlovassi, GR-83 200 Samos, Greece}
\address{School of Data Science and Artificial Inteligence,
Dongbei University of Finance and Economy, Dalian, Liaoning, China}
\email{konstant@aegean.gr} \email{sasd24009@sas.aegean.gr} \email{mengyuan\_prob@126.com}

\date{{\small \today}}

\begin{abstract}
In this paper we study local uniform, with respect to time, asymptotic expressions for the asymptotic behavior of the entrance probability of discounted aggregate claims to some rare sets, in a multivariate risk model with arbitrarily dependent insurance and financial risks. Our model is based on a multivariate version of a model, introduced by Guo (2022), and we consider that the logarithmic return process of the insurer's investment portfolio is described by a jump-diffusion process. The dependence between the insurance and financial risks is implied by the dependence of the claim-vectors with the jumps of returns, and is arbitrary under some distributional conditions on the claim-vectors and the 'discounted' claim-vectors. In opposite to previous papers on this topic, except the multi-dimensional extension, we consider that the model is driven by a common counting process, that is not necessarily renewal,  the claim vectors are interdependent, and their distribution is not restricted to the class of (multivariate) regularly varying distributions. Under the condition that the claim vectors, and the 'discounted' claim-vectors follows distributions from the class of multivariate consistently varying and positively decreasing distributions, and under some moment conditions on the jumps and the counting process, our main result shows the presence of multivariate linear single big jump principle of the discounted aggregate claims in this not necessarily L\'{e}vy-Renewal environment. After restriction of the distributions to multivariate regular variation, we obtain more explicit expressions, under a slightly weaker moment condition on the jumps. We provide a corollary, in which the conditions of the main result are satisfied under a weak dependence structure and we find a more explicit asymptotic expression, using the technique of solution of the dependence. Finally, under the multivariate regular variation frame, we give two  corollaries with weak and strong dependence structure, that depict directly the effect of the dependence between the insurance and financial risks on the insurer's solvency. Some examples with nested copulas included in order to stress the generality and the validity of our results.
\end{abstract}

\maketitle
\textit{Keywords: Multivariate risk model; Systemic factors; Non L\'{e}vy-Renewal environment; Interdependence; Multivariate linear single big jump}
\vspace{3mm}

\textit{Mathematics Subject Classification}: Primary 62P05 ;\quad Secondary 60G70.


\vspace{3mm}
\begin{center}
\textbf{ Dedicated to the memory of Professor \\{Alexander Alekseevich Borovkov (1931 - 2026)}}
\end{center}

\section{Introduction} \label{sec.KPT29U.1}

\subsection{Model description} \label{subsec.KPT29U.1.1}

In  modern actuarial practice the insurance companies invest their surplus into risk-free and risky assets for rising their solvency and competitiveness. Furthermore, for sake of risk diversification, they keep multiple portfolios, see for example \cite[p. 122]{aig:annual:report:2021}. For these reasons, many researchers focused recently on the asymptotic behavior of the ruin probability (or, related quantities) on multivariate risk models, in the presence of investments, see for example \cite{li:2016}, \cite{chen:li:cheng:2023}, \cite{yang:su:2023}, \cite{konstantinides:passalidis:2025j}, \cite{yuan:lu:fu:2025}, among others. These papers examine the case when the insurance risks, described by the claims, and the financial risks, represented by the return processes of the investment portfolio, are independent each other. Such an assumption is due to its mathematical tractability, and not in its practical verification. 
The world-scale crises, and the pandemics, indicate that such an assumption is too strict for the actuarial applications.

Very recently, some papers consider the study of continuous-time risk models in presence of dependence between insurance and financial risks, see in subsection \ref{subsec.KPT29U.1.2} for more discussions. A risk model, that seems to describe well periods of major financial instability, was suggested by \cite{guo:2022}, and in our paper we study a kind of multivariate analogue of that model. Concretely, we assume that the systemic factors arrive at time moments $\{\tau_{i},\,i \in \bbn\}$, with $\tau_0 =0$, and constitute a counting process 
\beam \label{eq.KPY.1.1}
N(t) = \sup \{i \in \bbn\;:\; \tau_i \leq t\}\,,
\eeam
for $t \geq 0$, with $\sup \emptyset =0$ conventionally. For $\{N(t)\,,\; t\geq 0\}$, we suppose that for any fixed $t\geq 0$, it holds $\lambda(t):=\E[N(t)]=\sum_{i=1}^{\infty} \PP(\tau_i \leq t) < \infty$. The systemic factors cause immediately, in the insurer's financial portfolio, jumps of sizes $\{Y_i\,,\; i \in \bbn\}$, that represent real valued random variables. We assume that the logarithmic return process of the investment portfolio is represented by the following jump diffusion process
\beam \label{eq.KPY.1.2}
\xi(t) =(r\,t + \sigma\,W(t)) + \sum_{i=1}^{N(t)} Y_i =: B(t) + J(t)\,,
\eeam 
for $t \geq 0$, where $r \in \bbr$ denotes the log-return rate, $\sigma > 0$, represents the volatility factor, while the $\{W(t)\,,\;t\geq 0\}$ is a Wiener process. According to \eqref{eq.KPY.1.2}, the process $\{B(t)\,,\;t\geq 0\}$ corresponds to a Brownian motion, with drift $r \in \bbr$, and volatility $\sigma >0$, hence its Laplace exponent at some $q \in \bbr$, is provided by the relation
\beao
\phi_B(s) := \log \E\left[ e^{-s\,B(1)}\right] = -r\,s + \dfrac{\sigma^2\,s^2}2\,,
\eeao
which is finite for any $s \in \bbr$. Hence by  \cite[Prop. 3.14]{cont:tankov:2004} we obtain
\beam \label{eq.KPY.1.3}
\E\left[ e^{-s\,B(t)}\right] =e^{t\,\phi_B(s)} < \infty\,,
\eeam 
for $t \geq 0$. Now, we adopt that the claim-vectors $\{{\bf X}^{(i)}\,,\;i \in\bbn\}$, arrive with some random delay $\theta \geq 0$ almost surely, from the systemic factors, hence their arrival times $\{\tau_i + \theta\,,\;i \in \bbn\}$ form the counting process
\beam \label{eq.KPY.1.4}
N_{\theta}(t):= \sup\{i \in \bbn\;:\; \tau_i + \theta \leq t\}\,,
\eeam
From relations \eqref{eq.KPY.1.1} and \eqref{eq.KPY.1.4} we see that for any fixed $t \geq 0$, the following inequality holds almost surely 
\beao
N_{\theta}(t) \leq N(t)\,.
\eeao
Consequently, we find $\lambda_{\theta}(t) :=\E[N_{\theta}(t)] \leq \lambda(t) < \infty $, for any fixed $t \geq 0$. To avoid trivialities, as the non-existence of renewal epochs at time $\tau_i + \theta$, we assume that each claim-vector ${\bf X}^{(i)}=( X_1^{(i)},\,\ldots,\, X_d^{(i)})^{\top}$, with $d\in \bbn$, (where ${\bf y}^{\top}$ denotes the transpose of vector ${\bf y}$), can has zero components, but not all the components zero. Additionally, we restrict our study at the time interval
\beam \label{eq.KPY.1.5}
\Lambda:= \{t \geq 0\;:\; \PP(\tau_1 + \theta \leq t) > 0\}\,,
\eeam
because, if $t \notin \Lambda$, then the probabilities of relation \eqref{eq.KPY.1.7} are equal to zero. For any fixed time $T \in \Lambda$, we denote $\Lambda_T := \Lambda \cap [0,\,T]$. So, the insurer's discounted aggregate claims at the moment $t \geq 0$ are given by the relation
\beam \label{eq.KPY.1.6}
&&{\bf D}(t)= \sum_{i=1}^{N_{\theta}(t)} {\bf X}^{(i)} \, e^{-\xi(\tau_i + \theta)} \\[2mm] \notag
&&= \left( \sum_{i=1}^{\infty} X_1^{(i)}\,e^{-\xi(\tau_i + \theta)}\,{\bf 1}_{\{\tau_i + \theta \leq t\}},\,\ldots,\,\sum_{i=1}^{\infty} X_d^{(i)}\,e^{-\xi(\tau_i + \theta)}\,{\bf 1}_{\{\tau_i + \theta \leq t\}} \right) ^{\top}\,,
\eeam
and we are interested in asymptotic estimation of the probabilities of the form
\beam \label{eq.KPY.1.7}
\PP({\bf D}(t) \in x\,A) \,,
\eeam
as $\xto$, where $A$ is a set from a general family of convex sets (see, below family $\mathscr{R}$ in subsection \ref{subsec.KPT29U.2.2}). The set $A$, can take several forms, that are useful in actuarial applications (see, relations \eqref{eq.KPY.2.4} and  \eqref{eq.KPY.2.5} below). 

The aim of this paper is the approximation of the probability at \eqref{eq.KPY.1.7}, with local uniformity with respect to time, namely uniformly over any set $\Lambda_T$, under weaker conditions, than that of previous papers, see in subsection \ref{subsec.KPT29U.1.2}. The main source of dependence between the insurance and financial risks, is depicted by the dependence among the pairs from $\{({\bf X}^{(i)},\,Y_i)\,,\;i \in \bbn\}$, while a secondary source of dependence between the two risk comes by the counting process $\{N(t)\,,\; t\geq 0\}$. The next two assumptions are valid through the whole paper.

\begin{assumption} \label{ass.KPY.1.1}
We suppose that the claims $\{{\bf X}^{(i)}\,,\;i \in\bbn\}$, are represented by non-negative, identically distributed random vectors, with common distribution $F$. The $\{Y_i\,,\;i \in \bbn \}$ are real-valued,  independent, but not necessarily identically distributed, random variables, and for any $q \geq 0$ we denote
\beao
M_q^{(i)} := \E\left[\left( e^{-Y_i}\right)^q \right]\,,
\eeao
for any $i \in \bbn$.  
\end{assumption}

\begin{assumption} \label{ass.KPY.1.2}
We suppose that the $\{N(t)\,,\; t\geq 0\}$, $\theta$, $\{W(t)\,,\; t\geq 0\}$, $\{({\bf X}^{(i)},\,Y_i)\,,\;i \in \bbn\}$, are mutually independent.   
\end{assumption}

These assumptions above, in one-dimensional set up, are used in \cite{guo:2022}, \cite{yang:fan:yuen:2023}, but as we shall see later in subsection \ref{subsec.KPT29U.1.2}, here we relax several of their assumptions, even in one-dimensional subcase, $d=1$.

Hereafter we organize the paper as follows. In subsection \ref{subsec.KPT29U.1.2} we provide a short overview of the recent literature on risk models with dependent insurance and financial risks, to define the relative contribution of this paper. In Section \ref{sec.KPT29U.2}, we present some introductory concepts, that concern the multivariate heavy-tailed distributions, which will be used in the main results, as also the dependence structure between claim vectors.

 In subsection \ref{subsec.KPT29U.3.1}, we present some main assumptions, that we need for Theorem \ref{th.KPY.3.1}, which refer to moment conditions on financial risks (concretely, on $M_q^{(i)}$), as also moment conditions on counting process $\{N(t)\,,\; t\geq 0\}$. We give some general examples of non-renewal processes, that satisfy this assumption. 

In subsection \ref{subsec.KPT29U.3.2} we provide our main result. In the case of restriction to distributions $F$ and $G_i$, with $i \in \bbn$, in the class of multivariate regular variation, $MRV$, we derive a more explicit asymptotic relation, than that of the main result, under some slightly weaker moment condition on $M_q^{(i)}$.

In Section \ref{sec.KPT29U.4} we present the proofs of the results, together with several necessary lemmas. 

Finally, in Section \ref{sec.KPT29U.5}, we present certain general examples and corollaries, in which the assumptions of our main results are satisfied. 
Among them, we provide a more direct asymptotic expression to our main theorem using the technique of the solution of dependence, while in the case of $MRV$ our two 
examples under weak and strong dependence structure bring light to the essential effect of the dependence of insurance and financial risks to insurer's solvency.

\subsection{Brief literature review} \label{subsec.KPT29U.1.2}

Recently, some papers have focused their attention on the study of continuous time risk models with dependence among the insurance and financial risks. The \cite{konstantinides:passalidis:2025o} was focused on the study of probability \eqref{eq.KPY.1.7} over infinite time horizon, with $t=\infty$. That risk model is more general than that of \eqref{eq.KPY.1.6}, with the assumption that the counting process is renewal and the logarithmic return process of the investment portfolio is a L\'{e}vy process, under several general dependencies and multivariate heavy-tailed distributions for the claim vectors. Further, \cite{konstantinides:passalidis:2026q} examined a similar model, under weak dependencies and smaller distribution classes in a non L\'{e}vy-Renewal framework. In spite of the fact that the assumptions are general enough, and appear as extensions of the Kesten - Goldie theorems, the case of finite horizon $t < \infty$, seems impossible to be handled by these techniques.

\cite{guo:2022} introduced the one-dimensional, $d=1$, model of \eqref{eq.KPY.1.6} for the study of ruin probability. Under Assumptions \ref{ass.KPY.1.1} and \ref{ass.KPY.1.2}, with $\{N(t)\,,\; t\geq 0\}$ to be a homogeneous Poisson process, the $\{Y_i\,,\;i \in \bbn\}$ independent and identically distributed (i.i.d) random variables, the $\{{\bf X}^{(i)}\,,\;i \in\bbn\}$ i.i.d. random vectors following regularly varying distribution with index $\alpha >0$, while $\widehat{\bf X}^{(i)}:={\bf X}^{(i)}\,e^{-Y_i} \stackrel{d}{\sim} G_i = G$, where $\stackrel{d}{\sim}$ indicates the distribution of the random vector $\widehat{\bf X}^{(i)}$, represents a regularly varying distribution, with index $\alpha >0$, and it provides uniform asymptotic estimations for the ruin probability for all the time horizons, using some conditions on $\phi_B(\cdot)$, \cite{yang:fan:yuen:2023} extends this study to the case when  $\{N(t)\,,\; t\geq 0\}$ is (homogeneous) renewal, the $G$ represents a regularly varying distribution with index not necessarily equal to $\alpha>0$, under some assumptions on $M_q^{(i)} = M_q $. These results hold local uniformly with respect to time. See also \cite{cheng:wang:2025}, \cite{xu:peng:zou:2025}, for relative papers on this topic.

In this paper we aim to asymptotic behavior of the discounted aggregate claims, namely of probability in \eqref{eq.KPY.1.7}, in the multivariate frame, and also along the following three extension directions
\begin{enumerate}
\item
The jumps $\{Y_i\,,\;i \in \bbn\}$ are not necessarily identically distributed.
\item
The $\{{\bf X}^{(i)}\,,\;i \in\bbn\}$ are not necessarily independent, but we use a general weak dependence structure, see $QAI_A$ in subsection \ref{subsec.KPT29U.2.3}. Further, the distribution $F$ and the distributions $G_i$, are not restricted only to multivariate regular variation, but we use a wider distribution class $(\mathcal{C}\cap \mathcal{P_D})_A$.
\item
Due to the fact that the model of \eqref{eq.KPY.1.6}, describes cases of high financial instability, the assumption that the $\{N(t)\,,\; t\geq 0\}$ is renewal seems too restrictive.
\end{enumerate}

The distribution class $(\mathcal{C}\cap \mathcal{P_D})_A$ can contain inhomogeneous risks among the insurance risks, which is not valid for the case of multivariate regular variation, where the marginal tails of the insurance risks should be mutually weakly equivalent (otherwise the Radon measure can be zero), something too restrictive for actuarial practice.

Here we suppose a general moment condition on $\{N(t)\,,\; t\geq 0\}$, that is satisfied by various counting processes, including all the renewal counting processes.

\section{Preliminaries} \label{sec.KPT29U.2}

In this Section we present some preliminary concepts related with the multivariate distributions with heavy tails, as also the dependence structures that we use among the claim vectors, after some necessary conventions.

\subsection{Notation} \label{subsec.KPT29U.2.1}

All the vectors are of dimension $d \in \bbn$, and denoted by bold script. For two vectors ${\bf x}$ and ${\bf y}$ we define their sum ${\bf x}+{\bf y}=(x_1+y_1,\,\ldots,\,x_d + y_d)^{\top}$ and for some constant $\lambda$, we define the scalar product $\lambda\,{\bf y}=(\lambda\,y_1,\,\ldots,\,\lambda\,y_d)^{\top}$, while by $|{\bf x}|$ we denote any norm in $\bbr_+^d$. For the real numbers $\lambda_1,\,\ldots,\,\lambda_n$, we define their maximum  $\bigvee_{i=1}^n \lambda_i :=\max\{\lambda_1,\,\ldots,\,\lambda_n\}$ and their minimum $\bigwedge_{i=1}^n \lambda_i :=\min\{\lambda_1,\,\ldots,\,\lambda_n\}$. For any set $\bbb$, we write $\bbb^c$ for the complement set, $\partial \bbb$ for the border, $\overline{\bbb}$ for the close hull, while with ${\bf 1}_{\bbb}$ we denote the indicator function of this set. Further, a set $\bbb$ is called increasing if for any ${\bf x} \in \bbb$ and ${\bf y} \in \bbr_+^d :=[0,\,\infty)^d$ it holds ${\bf x}+{\bf y} \in  \bbb$. We denote with ${\bf 0}=(0,\,\ldots,\,0)^{\top}$ the origin of the axes in $\bbr^d :=(-\infty,\,\infty)^d$.

Hereafter, all the limit relations hold as $\xto$, except is mentioned differently. For two positive $d$-dimensional functions ${\bf g}(\cdot)$ and ${\bf h}(\cdot)$ and for some set $\bbb \in \bbr^d$, with ${\bf 0} \notin \overline{\bbb}$, we denote ${\bf h}(x\,\bbb) \sim c\,{\bf g}(x\,\bbb)$, for some $c \in (0,\,\infty)$,  ${\bf h}(x\,\bbb) =O[{\bf g}(x\,\bbb)]$ and  ${\bf h}(x\,\bbb) =o[{\bf g}(x\,\bbb)]$, if it holds
\beao
\lim \dfrac{{\bf h}(x\,\bbb)}{{\bf g}(x\,\bbb)}=c\,,\quad \limsup \dfrac{{\bf h}(x\,\bbb)}{{\bf g}(x\,\bbb)}< \infty\,,\quad \lim \dfrac{{\bf h}(x\,\bbb)}{{\bf g}(x\,\bbb)}=0\,,
\eeao
respectively. We also write ${\bf h}(x\,\bbb) \asymp {\bf g}(x\,\bbb)$ if hold simultaneously ${\bf h}(x\,\bbb) =O[ {\bf g}(x\,\bbb)]$ and ${\bf g}(x\,\bbb) =O[ {\bf h}(x\,\bbb)]$. Further, for two $(d+1)$-variate positive functions ${\bf g}^*(\cdot)$ and ${\bf h}^*(\cdot)$, we say that ${\bf h}^*(x\,\bbb\,;\,z) \sim c\,{\bf g}^*(x\,\bbb\,;\,z)$, with $c \in (0,\,\infty)$, uniformly for $z \in E$, where $E \neq \emptyset$, if it holds
\beao
\lim \sup_{z\in E} \left|\dfrac{{\bf h}^*(x\,\bbb\,;\,z)}{c\,{\bf g}^*(x\,\bbb\,;\,z)}-1\right| =0\,.
\eeao
We declare that a random variable (or, vector) $Z$ follows the distribution $V$, as $Z \stackrel{d}{\sim} V$, while their support is denoted by $S(Z)$. If $V$ is one-dimensional, then we denote its tail by $\bV(x):= 1- V(x)$, for any $x \in \bbr$.

\subsection{Multivariate heavy-tailed distributions} \label{subsec.KPT29U.2.2}

For sake of compactness, all the distributions below are supported on $\bbr_{+}^d$. We note that in this subsection, all the one-dimensional distributions, namely $V$ and $F_A$, have infinite right endpoint, that means $\bF_A(x) \wedge \bV(x) >0$, for $x \in \bbr$.

For introduction of multivariate distribution classes, we need the set family
\beam \label{eq.KPY.2.1}
\mathscr{R} := \{ A \subsetneq \bbr^d\,:\, A \;{\text open, increasing}\,,\;A^c {\text convex}\,,\; {\bf 0} \notin \overline{A} \}\,,
\eeam
introduced in \cite{samorodnitsky:sun:2016}, for the definition of the multivariate subexponentiality. It was also proved that for any fixed $A \in \bbr$, if ${\bf X} \stackrel{d}{\sim} F$, with $F$ defined on $\bbr_+^d$, then for the random variable
\beam \label{eq.KPY.2.2}
X_A := \sup \{ u\;:\;{\bf X} \in u\,A\} \stackrel{d}{\sim} F_A\,,
\eeam
the $F_A$ represents a proper distribution, whose tail can be given by the relation
\beam \label{eq.KPY.2.3}
\bF_A(x) = \PP\left( \sup_{{\bf p} \in I_A}  {\bf p}^{\top}\,{\bf X}> x \right) = \PP({\bf X} \in x\,A)\,,
\eeam
for $x \geq 0$, where $I_A \subset \bbr^d$ represents an index set whose existence was proved for all $A \in \mathscr{R}$, see \cite[Lem.4.3 (c), Lem. 4.5]{samorodnitsky:sun:2016}, for the argumentation of the above.

Thus, by \eqref{eq.KPY.2.2} and \eqref{eq.KPY.2.3} was defined the multivariate subexponentiality on $A \in \mathscr{R}$ when $F_A$ is subexponential. Such a definition is closely related with the ruin probability, which in multivariate risk models can be defined by several ways. Two sets from family $\mathscr{R}$, that are useful in actuarial practice, are the set
\beam \label{eq.KPY.2.4}
A_1 = \{ {\bf y}\;:\;y_j > b_j\,,\; \exists j=1,\,\ldots,\,d \}\,,
\eeam
where $b_j >0$, for $j=1,\,\ldots,\,d$, and the set
\beam \label{eq.KPY.2.5}
A_2 = \left\{ {\bf y}\;:\;\sum_{j=1}^d l_j\,y_j > b \right\}\,,
\eeam
where $b>0$, $l_1,\,\ldots,\,l_d \geq 0$ and $l_1+\,\cdots+\,l_d=1 $. On these sets the probabilities \eqref{eq.KPY.1.7} play crucial role for the insurer's solvency. Also for $d=1$, we obtain $A_1=A_2 =(b,\,\infty)$, and hence \eqref{eq.KPY.1.7} is reduced to the distribution tail of the discounted aggregate claims. In this paper, we use a small subclass of multivariate subexponential distribution, to describe the claim vector distribution, more precisely we employ class $(\mathcal{C}\cap\mathcal{P_D})_A$.

In \cite{konstantinides:passalidis:2024g}, following the path of \cite{samorodnitsky:sun:2016}, introduced the following two classes. Let $A \in \mathscr{R}$ some fixed set. We say that the ${\bf X} \stackrel{d}{\sim} F$, belongs to the class of multivariate dominated varying distributions on $A$, symbolically $F \in \mathcal{D}_A$, if $F_A \in \mathcal{D}$, that means for any (or, equivalently, for some) $v\in (0,\,1)$ it holds
\beao
 \limsup \dfrac{\bF_A(v\,x)}{\bF_A(x)} < \infty\,.
\eeao  
We say that the ${\bf X} \stackrel{d}{\sim} F$, belongs to the class of consistently varying distributions on $A$, symbolically $F \in \mathcal{C}_A$, if $F_A \in \mathcal{C}$, that means it holds
\beao
\lim_{b\downarrow 1} \liminf \dfrac{\bF_A(b\,x)}{\bF_A(x)}=1\,.
\eeao
Further, \cite{konstantinides:passalidis:2025h} introduced the class $(\mathcal{P_D})_A$ of positively decreasing distributions on $A$. Concretely, we say that $F \in (\mathcal{P_D})_A$ if $F_A \in \mathcal{P_D}$, that means for any (or, equivalently, for some) $b>1$ it holds
\beao
\limsup \dfrac{\bF_A(b\,x)}{\bF_A(x)} < 1\,.
\eeao
Class $\mathcal{P_D}_A$ contains distributions either with heavy tails or with light tails. We say that $F \in (\mathcal{C}\cap \mathcal{P_D})_A$, if $F_A \in \mathcal{C}\cap \mathcal{P_D}$. Due to the fact that class $\mathcal{P_D}$ is a general enough class, the $\mathcal{C}_A \setminus (\mathcal{C}\cap \mathcal{P_D})_A$ is a too limited one, hence it does not contain distributions that are useful in actuarial practice. However, as we shall see later in Section \ref{sec.KPT29U.4}, this small restriction from class $\mathcal{C}_A$ to class $(\mathcal{C}\cap \mathcal{P_D})_A$ is necessary, to guarantee the non-defectiveness in the distribution of randomly weighted sums with infinite summands, that are needed for the argumentation of the main results. For further discussions about properties of class $\mathcal{P_D}$ we refer the reader to \cite{tang:2006a}, \cite{bardoutsos:konstantinides:2011} and \cite{konstantinides:passalidis:2024d}.

For all previously mentioned distribution classes, we denote their multivariate analogue on whole $\mathscr{R}$ as $\mathcal{B}_{\mathscr{R}} := \bigcap_{A \in \mathscr{R}} \mathcal{B}_A$, where $\mathcal{B} \in \{\mathcal{C},\,\mathcal{D},\,\mathcal{P_D},\,\mathcal{C}\cap\mathcal{P_D}\}$. For examples from these classes, see, \cite[Sec. 4]{konstantinides:liu:passalidis:2025} and Example \ref{exam.KPY.2.1} below.

In order to define the next class, we need first to define the one-dimensional regularly varying distributions. A distribution $V$ is called regularly varying with index $\alpha \in (0,\,\infty)$, symbolically $V \in \mathcal{R}_{-\alpha}$, if for any $y>0$ it holds
\beao
\lim \dfrac{\bV(y\,x)}{\bV(x)} =y^{-\alpha}\,.
\eeao
We say that ${\bf X} \stackrel{d}{\sim} F$ belongs to the class of multivariate regularly varying distributions, if there exists $V \in \mathcal{R}_{-\alpha}$, with  $\alpha \in (0,\,\infty)$, and a Radon measure $\mu(\cdot)$, non-degenerate to zero, such that it holds
\beam \label{eq.KPY.2.6}
\lim \dfrac{\PP({\bf X} \in x\,\bbb)}{\bV(x)} = \mu (\bbb)\,,
\eeam
for any Borel set $\bbb \in [0,\,\infty]^d$, with ${\bf 0} \notin \overline{\bbb}$, and $\mu(\partial \bbb) =0$. In this case we write $F \in MRV(\alpha,\,V,\,\mu)$. By the proof of \cite[Prop. 4.14]{samorodnitsky:sun:2016}, we have that if $F \in MRV(\alpha,\,V,\,\mu)$, then for any $A \in \mathscr{R}$ it holds $\mu(A) \in (0,\,\infty)$ and therefore from relation \eqref{eq.KPY.2.6} and by the closure property of the regular variation with respect to strong tail equivalence (see, \cite[Prop. 3.3 (i)]{leipus:siaulys:konstantinides:2023}), we obtain $F_A \in \mathcal{R}_{-\alpha}$. Hence, from the one-dimensional inclusions of the distribution classes with heavy tails (see, \cite[Ch. 2]{leipus:siaulys:konstantinides:2023}), we find immediately that
\beam \label{eq.KPY.2.7}
MRV \subsetneq (\mathcal{C}\cap \mathcal{P_D})_{\mathscr{R}} \subsetneq \mathcal{C}_{\mathscr{R}} \subsetneq \mathcal{D}_{\mathscr{R}}\,,
\eeam
where $MRV$ denotes the class of all multivariate regularly varying distributions. Inclusion \eqref{eq.KPY.2.7} holds obviously also for the classes $\mathcal{B}_A$ for any $A \in \mathscr{R}$, instead of $\mathcal{B}_\mathscr{R}$, where $\mathcal{B} \in \{\mathcal{C},\,\mathcal{D},\,\mathcal{C}\cap\mathcal{P_D}\}$. We refer the reader to \cite[Exam. 4.3, 4.4]{konstantinides:liu:passalidis:2025}, which show that the first inclusion in \eqref{eq.KPY.2.7} is not trivial. That examples focused mainly in the case of set $A_2$ of \eqref{eq.KPY.2.5}. Next, we provide an example which focus on the set $A_1$ of \eqref{eq.KPY.2.4}.

\bexam \label{exam.KPY.2.1}
Let $A_1$ defined by \eqref{eq.KPY.2.4}, ${\bf X}=(X_1,\,\ldots,\,X_d)^{\top} \stackrel{d}{\sim} F$. We suppose that the marginals $X_i \stackrel{d}{\sim} F_i \in \mathcal{B} \in \{\mathcal{C},\,\mathcal{C}\cap\mathcal{P_D}\}$ and are quasi-asymptotically independent ($QAI$), that means it holds
\beao
\lim \dfrac{\PP(X_i > x\,,\;X_j > x)}{\bF_i(x) + \bF_j(x)} =0\,,
\eeao
for any $1\leq i \neq j \leq d$. Then $F \in \mathcal{B}_{A_1}$, with $\mathcal{B} \in \{\mathcal{C},\,\mathcal{C}\cap\mathcal{P_D}\}$, respectively.

\pr~
Since the $(X_1,\,\ldots,\,X_d)$ are $QAI$ and non-negative, we obtain
\beam \label{eq.KPY.2.8}
\bF_{A_1}(x) = \PP({\bf X} \in x\,A_1) = \PP\left(\bigvee_{i=1}^d \dfrac{X_i}{b_i} > x\right) \leq \PP\left(\sum_{i=1}^d \dfrac{X_i}{b_i} > x\right)  \sim \sum_{i=1}^d \PP\left( \dfrac{X_i}{b_i} > x\right) \,,
\eeam
where at the last step we used \cite[Th. 3.2]{chen:yuen:2009}. From the other hand side, by Bonferroni inequality and taking into account the $QAI$ property we find
\beam \label{eq.KPY.2.9}
&&\bF_{A_1}(x) = \PP\left(\bigvee_{i=1}^d \dfrac{X_i}{b_i} > x\right) \\[2mm] \notag
&& \geq \sum_{i=1}^d \PP\left( \dfrac{X_i}{b_i} > x\right) -  \sum_{1 \leq i< j \leq d} \PP\left( \dfrac{X_i}{b_i} > x\,,\;\dfrac{X_j}{b_j} > x\right) \sim \sum_{i=1}^d \PP\left( \dfrac{X_i}{b_i} > x\right) \,.
\eeam
If $X_i/b_i \stackrel{d}{\sim} F_i'$, for $i =1,\,\ldots,\,d$, then obviously $F_i' \in \mathcal{B} \in \{\mathcal{C},\,\mathcal{C}\cap\mathcal{P_D} \}$, for any $i =1,\,\ldots,\, d$. So, in combination with the asymptotic relation
\beam \label{eq.KPY.2.b}
\bF_{A_1}(x) \sim \sum_{i=1}^d \PP\left( \dfrac{X_i}{b_i} > x\right)=\sum_{i=1}^d \bF_i'(x)\,,
\eeam
that follows from \eqref{eq.KPY.2.8} and \eqref{eq.KPY.2.9}, it is easy to see that $F_{A_1} \in \mathcal{B} \in \{\mathcal{C},\,\mathcal{C}\cap\mathcal{P_D}\}$ and hence  $F \in \mathcal{B}_{A_1} \in \{\mathcal{C}_{A_1},\,(\mathcal{C}\cap\mathcal{P_D})_{A_1}\}$.
~\halmos
\eexam

We observe in the previous example that we do not need any assumption about 'weak equivalence' of tails for the marginal distributions, something that is necessary for the $MRV$ case, except of cases in the form $\mu([0,\,\infty],\,\ldots,\,(1,\,\infty],\,\ldots,\,[0,\,\infty])=0$, which indicates that some (or, more) components have asymptotically negligible tail in comparison with $V \in \mathcal{R}_{-\alpha}$. However, in theses cases $MRV$ is not helpful with respect of its use, since does not provide sharp asymptotic estimations. The approach by \cite{samorodnitsky:sun:2016} however seems to provide better approximation in such cases of inhomogeneous risks, as follows from Example \ref{exam.KPY.2.1}.  

We close this subsection presenting the Matuszewska indexes, that are helpful in the characterization of several distribution classes with heavy tails, and we shall need the moment conditions of jumps of financial risks, namely $M_q^{(i)}$. The lower and upper Matuszewska indexes of an one-dimensional distribution $V$, are respectively defined as 
\beao
J_V^- = -\lim_{b \to \infty} \dfrac{\log \overline{V^*}(b)}{\log b}\,, \quad J_V^+ = -\lim_{b \to \infty} \dfrac{\log \overline{V_*}(b)}{\log b}\,,
\eeao
where 
\beao
\overline{V^*}(b) = \limsup \dfrac{ \overline{V}(b\,x)}{\overline{V}(x)}\,, \quad \overline{V_*}(b) = \liminf \dfrac{\overline{V}(b\,x)}{\overline{V}(x)}\,.
\eeao
For all distributions $V$ it holds $0 \leq J_V^- \leq J_V^+ \leq \infty$. It is well known that $V \in \mathcal{D}$, if and only if $J_V^+ < \infty$, if $V \in \mathcal{R}_{-\alpha}$, with $\alpha \in (0,\,\infty)$, then $J_V^- = J_V^+ =\alpha$, while $V \in \mathcal{P_D}$ if and only if $J_V^->0$. See \cite[Sec. 2.1]{bingham:goldie:teugels:1987}, \cite[Subsec. 2.4]{leipus:siaulys:konstantinides:2023}, for more discusions on this indexes.

By \cite[Lem. 3.5]{tang:tsitsiashvili:2003} we obtain that if $J_V^+ < \infty$, then for any fixed $p>J_V^+$ it holds
\beam \label{eq.KPY.2.10}
x^{-p} = o[\bV(x)] \,.
\eeam

\subsection{Quasi asymptotic independence on $A$} \label{subsec.KPT29U.2.3}

In this subsection, we describe the dependence structure, that we employ among the claim vectors $\{{\bf X}^{(i)}\,,\;i \in \bbn\}$. For any $i \in \bbn$, we denote
\beam \label{eq.KPY.2.11}
X_A^{(i)} := \sup\{ u\;:\; {\bf X}^{(i)} \in u\,A\} \stackrel{d}{\sim} F_A^{(i)} \,.
\eeam
The following dependence structure, is a multivariate analogue of the $QAI$, introduced by \cite{chen:yuen:2009}.

\bde \label{def.KPY.2.1}
Let $A \in \mathscr{R}$ some fixed set. We say that $\{{\bf X}^{(i)}\,,\;i \in \bbn\}$ are quasi asymptotically independent on $A$, symbolically $QAI_A$, if the $\{X_A^{(i)}\,,\;i \in \bbn\}$ are $QAI$, namely for any $i,\,j \in \bbn$ with $i\neq j$ it holds
\beam \label{eq.KPY.2.12}
\lim \dfrac{\PP(X_A^{(i)} > x\,,\;X_A^{(j)} > x)}{\PP(X_A^{(i)} > x)+\PP(X_A^{(j)} > x)} = 0\,.
\eeam
\ede

\bre \label{rem.KPY.2.1}
In the framework of one-dimensional set up, $QAI$ has been studied in several papers, and is satisfied by many commonly used copulas, see for example \cite{chen:yuen:2009}, \cite{li:2013}, \cite{cheng:2014} among others. From Definition \ref{def.KPY.2.1}, we can see that if $\{{\bf X}^{(i)}\,,\;i \in \bbn\}$ are independent, then they are also $QAI_A$, for any $A \in \mathscr{R}$. In \cite[Sec. 2.2]{chen:konstantinides:passalidis:2025}, we obtain two examples of $QAI_A$, where $A= A_1'$.

In the following example we show sufficient conditions for the $QAI_{A_1'}$, in combination with the distributions of the claim vectors in $(\mathcal{C}\cap \mathcal{P_D})_{A_1'}$. So, this example contains sufficient conditions for the sequence $\{{\bf X}^{(i)}\,,\;i \in \bbn\}$, that appear in the main result, see Theorem \ref{th.KPY.3.1} below.

For sake of simplicity, we use the set 
\beam \label{eq.KPY.2.13}
A_1' = \{(x,\,y)\;:\; x\vee y > 1\}\,,
\eeam  
but the result holds through similar methodology for any $b>0$, instead of $b=1$, as in \eqref{eq.KPY.2.13}.
\ere

\bexam \label{exam.KPY.2.2}
Let ${\bf X}=(X_1,\,X_2)$ and ${\bf Z}=(Z_1,\,Z_2)$, two non-negative random vectors, with marginal distributions from class $\mathcal{B} \in \{\mathcal{C},\,\mathcal{C}\cap\mathcal{P_D}\}$. We suppose that the $X_1,\,X_2,\,Z_1,\,Z_2$, are widely upper orthant dependent, symbolically $WUOD$, see for examples of this dependence in \cite{wang:wang:gao:2013}, \cite{wang:cheng:2011}. Namely, for the $Y_1,\,Y_2,\,Y_3,\,Y_4$ with $Y_i \in \{X_1,\,X_2,\,Z_1,\,Z_2\}$, for $i =1,\,\ldots,\,4$, with $Y_i \neq Y_j$, for any $1\leq i \neq j \leq 4$, there exists a sequence of positive constants $\{g_U(n)\,,\;n =1,\,\ldots,\,4\}$ such that it holds
\beam \label{eq.KPY.2.14}
\PP\left(\bigcap_{i=1}^n \{Y_i > x_i\} \right) \leq g_U(n) \prod_{i=1}^n \PP(Y_i >x_i)\,,
\eeam 
for all $x_i \in \bbr$, with $i=1,\,\ldots,\,n$. Then the ${\bf X}$ and  ${\bf Z}$ follow distributions from the class $\mathcal{B}_{A_1'}$, with $\mathcal{B} \in \{\mathcal{C},\,\mathcal{C}\cap\mathcal{P_D}\}$ and are $QAI_{A_1'}$.

\pr~
Firstly, it holds that $WUOD \subsetneq QAI$. Indeed, since by \eqref{eq.KPY.2.14}, for any $1\leq i \neq j \leq 4$ we obtain
\beam \label{eq.KPY.2.15} \notag
&&\lim \dfrac{\PP\left(Y_i > x\,,\;Y_j > x \right)}{\PP\left(Y_j > x \right)} \leq g_U(2) \lim \dfrac{\PP\left(Y_i > x \right)\,\PP\left(Y_j > x \right)}{\PP\left(Y_j > x \right)} \\[2mm]
&&=g_U(2) \lim \PP\left(Y_i > x \right)=0\,,
\eeam 
and from \eqref{eq.KPY.2.15} we take immediately that $QAI$ is valid. Hence, from Example \ref{exam.KPY.2.1}, we find that the distributions of ${\bf X}$ and ${\bf Z}$ belong to $\mathcal{B}_{A_1'}$, with $\mathcal{B} \in \{\mathcal{C},\,\mathcal{C}\cap\mathcal{P_D}\}$. Further, by \eqref{eq.KPY.2.b} we obtain
\beam \label{eq.KPY.2.16}
\PP\left(X_{A_1'}> x \right)\sim \sum_{i=1}^2 \PP\left(X_i > x \right)\,,\quad \PP\left(Z_{A_1'}> x \right)\sim \sum_{i=1}^2 \PP\left(Z_i > x \right)\,,
\eeam 
where only for this example we denote $X_{A_1'}:= \sup \{u\;:\; {\bf X} \in u\,A_1' \}$, $Z_{A_1'}:= \sup \{u\;:\; {\bf Z} \in u\,A_1' \}$. It remains to show that the ${\bf X}$ and ${\bf Z}$, are $QAI_{A_1'}$, that means to prove the limit
\beam \label{eq.KPY.2.17}
\lim \dfrac{\PP\left(X_{A_1'}> x\,,\;Z_{A_1'}> x \right)}{\PP\left(X_{A_1'}> x \right)+\PP\left(Z_{A_1'}> x \right)} =0\,,
\eeam
For any $x>0$ we obtain
\beam \label{eq.KPY.2.18}
&&\PP\left(X_{A_1'}> x\,,\;Z_{A_1'}> x \right)=\PP\left(\{X_1> x\} \cup \{X_2> x\}\,,\;\{Z_1> x\} \cup \{Z_2> x\}\right)\\[2mm] \notag
&&=\PP\left(\{X_1> x\} \cup \{X_2> x\}\right)+\PP\left(\{Z_1> x\} \cup \{Z_2> x\}\right)- \\[2mm] \notag
&&\PP\left(\{X_1> x\} \cup \{X_2> x\}\cup \{Z_1> x\} \cup \{Z_2> x\}\right)=: \PP(X_{A_1'}> x)+\PP(Z_{A_1'}> x)-P_3\,,
\eeam
For $P_3$ we find that
\beam \label{eq.KPY.2.19}
&& P_3=\sum_{i=1}^2 \PP\left(X_i> x\right) +\sum_{i=1}^2 \PP\left(Z_i> x\right)- \PP\left(X_1> x\,,\;X_2> x\right)\\[2mm] \notag
&& - \PP\left(X_1> x\,,\;Z_1> x\right) - \PP\left(X_2> x\,,\;Z_1> x\right) -\PP\left(Z_1> x\,,\;Z_2> x\right) \\[2mm] \notag
&& - \PP\left(X_2> x\,,\;Z_2> x\right) - \PP\left(X_1> x\,,\;Z_2> x\right) \\[2mm] \notag
&& + \PP\left(X_1> x\,,\;X_2> x\,,\;Z_1> x\right)+ \PP\left(X_1> x\,,\;X_2> x\,,\;Z_2> x\right) \\[2mm] \notag
&&+ \PP\left(X_1> x\,,\;Z_1> x,,\;Z_2> x\right)+ \PP\left(X_2> x\,,\;Z_1> x\,,\;Z_2> x\right)  \\[2mm] \notag
&&-  \PP\left(X_1> x\,,\;X_2> x\,,\;Z_1> x\,,\;Z_2> x\right)\\[2mm] \notag
&& \geq \sum_{i=1}^2 \PP\left(X_i> x\right) +\sum_{i=1}^2 \PP\left(Z_i> x\right)- \PP\left(X_1> x\,,\;X_2> x\right)-\PP\left(X_1> x\,,\;Z_1> x\right) \\[2mm] \notag
&& -\PP\left(X_2> x\,,\;Z_1> x\right)-  \PP\left(Z_1> x\,,\;Z_2> x\right) - \PP\left(X_2> x\,,\;Z_2> x\right) \\[2mm] \notag
&&- \PP\left(X_1> x\,,\;Z_2> x\right)- \PP\left(X_1> x\,,\;X_2> x\,,\;Z_1> x\,,\;Z_2> x\right) \,,
\eeam
Due to \eqref{eq.KPY.2.14}, similarly with \eqref{eq.KPY.2.15}, for any $\vep > 0$, we can find a sufficiently large $x_0 = x_0(\vep) > 0$, such that for any $x \geq x_0$ it holds
\beao
\dfrac{\PP\left(Y_i> x\,,\;Y_j> x\right)}{\PP\left(Y_i> x\right)} < \vep\,,
\eeao
for any $1 \leq i\neq j \leq 4$ and similarly
\beao
\dfrac{\PP\left(X_1> x\,,\;X_2> x\,,\;Z_1> x\,,\;Z_2> x\right)}{\PP\left(X_1> x\right)} < \vep\,.
\eeao
From \eqref{eq.KPY.2.19} and the last two relations is implied that for any $x \geq x_0$ it holds
\beam \label{eq.KPY.2.20} \notag
&& P_3 \geq \sum_{i=1}^2 \PP\left(X_i> x\right) +\sum_{j=1}^2 \PP\left(Z_j> x\right)- 7\,\vep\,\left(\sum_{i=1}^2 \PP\left(X_i> x\right) +\sum_{j=1}^2 \PP\left(Z_j> x\right) \right) \\[2mm]
&& =(1-7\,\vep) \left(\sum_{i=1}^2 \PP\left(X_i> x\right) +\sum_{j=1}^2 \PP\left(Z_j> x\right) \right) \,,
\eeam
From \eqref{eq.KPY.2.16}, \eqref{eq.KPY.2.18} and \eqref{eq.KPY.2.20}, we obtain that for any $x \geq x_0$ it holds
\beao
&&\PP\left(X_{A_1'}> x\,,\;Z_{A_1'}> x \right)   \\[2mm]
&&\leq \PP\left(X_{A_1'}> x \right)+ \PP\left(Z_{A_1'}> x\right)- (1-7\,\vep)\,(1-\vep)\,\left[ \PP\left(X_{A_1'}> x \right)+ \PP\left(Z_{A_1'}> x\right)\right] \\[2mm]
&&=[1- (1-8\,\vep +7\,\vep^2)]\,\left[ \PP\left(X_{A_1'}> x \right)+ \PP\left(Z_{A_1'}> x\right)\right]\\[2mm]
&&=[8\,\vep - 7\,\vep^2]\,\left[ \PP\left(X_{A_1'}> x \right)+ \PP\left(Z_{A_1'}> x\right)\right]\,.
\eeao 
From the last inequality, given the arbitrary choice of $\vep >0$, we conclude that \eqref{eq.KPY.2.17} is valid.
~\halmos
\eexam

\section{Main results} \label{sec.KPT29U.3}

In this section we present the main result, after the introduction of two basic assumptions. In case, where we restrict the distributions $F$, $G_i$, with $i \in \bbn$, into class $MRV$, are derived much more explicit asymptotic expressions.

\subsection{Moment conditions} \label{subsec.KPT29U.3.1}

\begin{assumption} \label{ass.KPY.3.1}
Let $A \in \mathscr{R}$, some fixed set. We suppose that $F \in (\mathcal{C}\cap \mathcal{P_D})_A$, $G_i \in (\mathcal{C}\cap \mathcal{P_D})_A$, for any $i \in \bbn$ and for $i\neq j$, it holds $G_i(x\,A) \asymp G_j(x\,A) $, with $i,\,j\; \in \bbn$. Further, let there exist $0< p_1 < J_{F_A}^- \leq J_{F_A}^+ < p_2 < \infty$ and $0< q_1 < J_{G_A^{(1)}}^- \leq J_{G_A^{(1)}}^+ < q_2 < \infty$, 
such that it holds
\beam \label{eq.KPY.3.1}  
\left( M_{p_1}^{(i)} \vee M_{p_2}^{(i)} \right) \bigvee \left( M_{q_1}^{(i)} \vee M_{q_2}^{(i)} \right) \leq 1\,,
\eeam 
for any $i \in \bbn$.   
\end{assumption}

\bre \label{rem.KPY.3.1}
We have to notice that the assumption $G_i(x\,A) \asymp G_j(x\,A)$, with $i,\,j\; \in \bbn$ and $i \neq j$, implies that $\bG_A^{(i)}(x) \asymp \bG_A^{(j)}(x)$, where $G_A^{(i)}$ represents the distribution of
\beam \label{eq.KPY.3.a}  
\widehat{X}_A^{(i)}:= \sup\{u\;:\;\widehat{\bf X}^{(i)} \in u\,A \}\,.
\eeam 
This means that for all $i \in\bbn$ it holds $J_{G_A^{(i)}}^+=J_{G_A^{(1)}}^+$ and $J_{G_A^{(i)}}^-=J_{G_A^{(1)}}^-$. 

Further, the relation \eqref{eq.KPY.3.1} of Assumption \ref{ass.KPY.3.1}  implies the fact that insurance risks dominate on the financial risks, in the sense that the $e^{\xi(\tau_{i-1})-\xi(\tau_i)}$ follows distribution with lighter tail than that of $F_A$. Indeed, for $p_2 > J_{F_A}^+$ we obtain 
\beam \label{eq.KPY.3.2}  
&&\E\left[ \left( e^{\xi(\tau_{i-1})-\xi(\tau_i)} \right)^{p_2} \right] = \E\left[ \left( e^{B(\tau_{i-1})-B(\tau_i)-Y_i} \right)^{p_2} \right] \\[2mm] \notag
&& =\E\left[ \left( e^{-Y_i} \right)^{p_2} \right]\,\E\left[ e^{-p_2[B(\tau_i-\tau_{i-1})]} \right] = M_{p_2}^{(i)}\,\E\left[ e^{-p_2[B(w_i)]} \right] \leq \E\left[ e^{w_i\,\phi_B(p_2)} \right] < \infty\,, 
\eeam
where $\{w_i:=\tau_i - \tau_{i-1}\,,\; i \in \bbn\}$ represents the inter-arrival times between successive claim arrivals. At the pre-last step we used relation \eqref{eq.KPY.1.3} 
and \eqref{eq.KPY.3.1}, while the last step follows by the fact that $\phi_B(\cdot)$ is finite, on any finite point. From \eqref{eq.KPY.3.2} we see that the same would hold and 
for $M_{p_2}^{(i)} < \infty$, namely $\E\left[ \left( e^{\xi(\tau_{i-1})-\xi(\tau_i)} \right)^{p_2} \right] < \infty$, however relation \eqref{eq.KPY.3.1} is necessary for the convergence of the infinite randomly weighted sums, that were used in the proof of Theorem \ref{th.KPY.3.1}. The domination of the insurance risks against the financial risks is an assumption that is used in most of the papers on risk theory, although such an assumption was not observed always in actuarial  practice, see \cite{li:tang:2015}, \cite{tang:yang:2019} and \cite{konstantinides:passalidis:2025o}, for some discussions on this opposite case.          
\ere

\bre \label{rem.KPY.3.2} 
We can see that when $\{Y_i\,,\;i \in \bbn\}$ are non-negative random variables, eventually degenerate to zero, then \eqref{eq.KPY.3.1} is valid immediately. In this case, the process $\{\xi(t)\,,\;t \geq 0\}$ represents the logarithmic returns in an 'almost' risk-free portfolio, in the sense that the fluctuations in returns seem rather small and the bad scenarios caused only by the Wiener process. In such kind of cases, the dependencies between ${\bf X}^{(i)}$ and $e^{-Y_i}$ are necessarily weak, and even negatively associated.  
\ere

For our model we need also the following mild moment assumption  for the counting process $\{N(t)\,,\; t\geq 0\}$. 

\begin{assumption} \label{ass.KPY.3.2}
For the $p_1,\,p_2,\,q_1,\,q_2$, from Assumption \ref{ass.KPY.3.1}, let denote $q^*=p_1\vee p_2\vee q_1\vee q_2$. We suppose that there exists some fixed $\gamma > 1$, such that for all fixed $T \in \Lambda$ it holds
\beam \label{eq.KPY.3.3}  
\E\left[ \left(N(t) \right)^{\gamma\,q^*+1} \right] < \infty\,,
\eeam 
for any $t \in \Lambda_T$.   
\end{assumption}

By \cite{stein:1946} we find that Assumption \ref{ass.KPY.3.2} is satisfied by all renewal processes, since any renewal process has analytical moment-generating function near to zero. We provide now two examples, where \eqref{eq.KPY.3.3} is satisfied for non-renewal processes. The first example refers to a quasi-renewal process, namely a process with identically distributed but not necessarily independent inter-arrival times $\{w_i\,,\; i \in\bbn\}$.

\bexam \label{exam.KPY.3.1} 
Let consider inter-arrival times $\{w_i\,,\; i \in\bbn\}$, that are identically distributed with expectation $\E[w_i] =1/\lambda$, for any $i \in \bbn$, and extended negatively dependent, symbolically $END$, namely there exists a constant $M>0$, such that  for any $n \in \bbn$ it holds
\beao
&&\PP\left( \bigcap_{i=1}^n \{w_i > x_i \}\right) \leq M\,\prod_{i=1}^n \PP(w_i > x_i)\,, \\[2mm] \notag 
&&\PP\left( \bigcap_{i=1}^n \{w_i \leq x_i \}\right) \leq M\,\prod_{i=1}^n \PP(w_i \leq x_i)\,,
\eeao
for all $x_i \in \bbr$, with $i=1,\,\ldots,\,n$, see \cite{liu:2009} and \cite{chen:chen:ng:2010}, for more discussions about this dependence. Then, for any fixed $T \in \Lambda$, \eqref{eq.KPY.3.3} is true for all $t \in\Lambda_T$.

\pr~
From \cite[Th. 1.2]{wang:cheng:2011}, for any fixed $\delta >0$, we obtain 
\beam \label{eq.KPY.3.4}  
\E\left[ \left(N(t) \right)^{\delta} \right] \sim (\lambda\,t)^{\delta}\,,
\eeam 
as $t \to \infty$. From \eqref{eq.KPY.3.4}, for any $\vep>0$, we can find some $t_0>0$, large enough, such that 
\beam \label{eq.KPY.3.5}  
\E\left[ \left(N(t_0) \right)^{\gamma\,q^*+1} \right] \leq (1+ \vep) (\lambda\,t_0)^{\gamma\,q^*+1} < \infty\,.
\eeam

If $t_0 \geq T$, then because of the inequality $\E\left[ \left(N(t) \right)^{\gamma\,q^*+1} \right] \leq \E\left[ \left(N(t_0) \right)^{\gamma\,q^*+1} \right]$, by \eqref{eq.KPY.3.5} follows immediately \eqref{eq.KPY.3.3} for any $t \in \Lambda_T$.

If $t_0 < T$, then \eqref{eq.KPY.3.5} still remains valid with $T$ instead of $t_0$. Hence because of the inequality $\E\left[ \left(N(t) \right)^{\gamma\,q^*+1} \right] \leq \E\left[ \left(N(T) \right)^{\gamma\,q^*+1} \right]$, relation \eqref{eq.KPY.3.3} holds for any $t \in \Lambda_T$.
~\halmos
\eexam

The second example is related with a wide spectrum of inhomogeneous renewal processes.

\bexam \label{exam.KPY.3.2} 
Let consider inter-arrival times $\{w_i\,,\; i \in\bbn\}$, that are independent random variables, that are uniformly integrable, namely it holds
\beao
\lim_{z \to \infty} \sup_{i \in \bbn} \E\left[ w_i \,{\bf 1}_{\{w_i \geq z\}}\right] = 0\,.
\eeao
If it  holds
\beao
\lim_{n \to \infty} \dfrac 1n \sum_{i=1}^n \E\left[ w_i \right] = \dfrac 1{\mu}\,,
\eeao
for some $\mu \in (0,\,\infty)$. Then, for any fixed $T \in \Lambda$, \eqref{eq.KPY.3.3} holds for any $t \in \Lambda_T$.

\pr~
From \cite[Cor. 2.2]{bernackaite:siaulys:2015}, for any fixed $\delta>0$, we obtain \eqref{eq.KPY.3.4}. Now, with quite similar way, as in Example \ref{exam.KPY.3.1}, we can find that \eqref{eq.KPY.3.3} holds for any $t \in \Lambda_T$.
~\halmos
\eexam

\bre \label{rem.KPY.3.3}
We note that from the line of proof of Examples \ref{exam.KPY.3.1} and \ref{exam.KPY.3.2}, recall \eqref{eq.KPY.3.4}, for any fixed $\delta >0$, we find that relation \eqref{eq.KPY.3.3} is satisfied for any $t \in \Lambda_T$, independently of the choice of $\gamma$ and $q^*$. 
\ere

\subsection{Main theorem } \label{subsec.KPT29U.3.2}

The following theorem, is the main result of this paper. For this we need the introduction of the random variables
\beam \label{eq.KPY.3.6}  \notag
&&\Theta_i(t)= e^{-B(\tau_i +\theta) -J[(\tau_i +\theta)-]} {\bf 1}_{\{\tau_i +\theta \leq t\}} = \exp\left\{-B(\tau_i +\theta) -\sum_{j=1}^{N[(\tau_i + \theta)-]} Y_j \right\} {\bf 1}_{\{\tau_i +\theta \leq t\}} \,,\\[2mm] 
&&\widetilde{\Theta}_i(t)= \exp\left\{-B(\tau_i +\theta) +Y_i -\sum_{j=1}^{N[(\tau_i + \theta)-]} Y_j \right\} {\bf 1}_{\{\tau_i +\theta \leq t\}} \,,
\eeam
for any $i \in \bbn$, with $t \in \Lambda$.

\bth \label{th.KPY.3.1}
Let consider the discounted aggregate claims from \eqref{eq.KPY.1.2}. Under Assumptions \ref{ass.KPY.1.1}, \ref{ass.KPY.1.2}, \ref{ass.KPY.3.1}, \ref{ass.KPY.3.2} and if the sequence $\{{\bf X}^{(i)}\,,\;i \in \bbn\}$ is $QAI_A$, it holds
\beam \label{eq.KPY.3.7}   \notag
&&\PP\left( {\bf D}(t) \in x\,A \right) \sim \sum_{i=1}^{\infty} \PP\left({\bf X}^{(i)}\,\Theta_i(t)\,{\bf 1}_{\{\theta =0\}}  \in x\,A \right) +\sum_{i=1}^{\infty} \PP\left(\widehat{\bf X}^{(i)}\,\widetilde{\Theta}_i(t)\,{\bf 1}_{\{\theta >0\}}  \in x\,A \right) \\[2mm]
&&= \sum_{i=1}^{\infty} \PP\left({\bf X}^{(i)}\,\Theta_i(t) \in x\,A \right) \,,
\eeam
uniformly for $t \in \Lambda_T$, for any fixed $T \in \Lambda$.
\ethe

\bre \label{rem.KPY.3.4}
Relation \eqref{eq.KPY.3.7} looks a bit complicated. However, it contains the products of independent random vectors ${\bf X}^{(i)}\,\Theta_i(t)\,{\bf 1}_{\{\theta =0\}}$ and $\widehat{\bf X}^{(i)}\,\widetilde{\Theta}_i(t)\,{\bf 1}_{\{\theta >0\}}$, while later we provide through Corollary \ref{cor.KPY.5.1} some cases, in which the $\widehat{\bf X}^{(i)}$ can be also approached via an independent product. Further, due to the proof of Theorem \ref{th.KPY.3.1}, relation \eqref{eq.KPY.3.7} contains a kind of uniformity with respect to summands, see Lemma \ref{lem.KPY.4.3} below, that helps to numerical approximations of the infinite sums of \eqref{eq.KPY.3.7}, by finite sums with large number of summands. This shows the strong  presence of the multivariate linear single big jump principle, for the discounted aggregate claims, although we are not limited on L\'{e}vy-Renewal environment.
\ere

If we restrict the distributions $F$, $G_i$, with $i \in \bbn$, into class $MRV$, then \eqref{eq.KPY.3.7} reaches to a more explicit form, in which the dependencies of each vector ${\bf X}^{(i)}$, $\widehat{\bf X}^{(i)}$, are described through the Radon measures of $MRV$, while the decay rate of probability \eqref{eq.KPY.1.7}, is determined by two one-dimensional regularly varying distributions. The following corollary presents this case of restriction into $MRV$, however with slightly weaker moment conditions, than that from \eqref{eq.KPY.3.1}.

\bco \label{cor.KPY.3.1} 
Let consider the discounted aggregate claims in relation \eqref{eq.KPY.1.2}. Under the conditions of Theorem \ref{th.KPY.3.1} with the restrictions $F \in MRV(\alpha,\,V,\,\mu)$, $G_i \in MRV(\beta,\,Q_i,\,\nu_i)$, with $i \in \bbn$, and $\alpha,\,\beta \in (0,\,\infty)$, and instead of  \eqref{eq.KPY.3.1}, we assume
\beam \label{eq.KPY.3.8}   
M_{\alpha}^{(i)} \vee M_{\beta}^{(i)} \leq 1\,, \qquad (M_{p_1}^{(i)} \vee M_{p_2}^{(i)} ) \bigvee (M_{q_1}^{(i)} \vee M_{q_2}^{(i)}) < \infty\,,
\eeam   
for any $i \in \bbn$, then it holds
\beam \label{eq.KPY.3.9}  
\PP\left( {\bf D}(t) \in x\,A \right) \sim \mu(A) \bV(x) \sum_{i=1}^{\infty} \E\left[\Theta_i^{\alpha}(t)\,{\bf 1}_{\{\theta =0\}} \right] +\sum_{i=1}^{\infty} \nu_i(A) \overline{Q}_i(x)\,\E\left[\widetilde{\Theta}_i^{\beta}(t)\,{\bf 1}_{\{\theta >0\}} \right] ,
\eeam
uniformly for $t \in \Lambda_T$, for any fixed $T \in \Lambda$.
\eco

\bre \label{rem.KPY.3.5}
As in Theorem \eqref{th.KPY.3.1}, here again the estimation in \eqref{eq.KPY.3.9} is in fact uniform with respect to $n \in \bbn$, see Lemma \ref{lem.KPY.4.4} below. We also note, that in that case the $\{Y_i\,,\; i \in \bbn\}$ are identically distributed, and hence $G_i = G \in MRV(\beta,\,Q,\,\nu) $, we obtain that relation \eqref{eq.KPY.3.9} is reduced to
\beam \label{eq.KPY.3.10}  
\PP\left( {\bf D}(t) \in x\,A \right) \sim \mu(A) \bV(x) \sum_{i=1}^{\infty} \E\left[\Theta_i^{\alpha}(t)\,{\bf 1}_{\{\theta =0\}} \right] +\nu(A) \overline{Q}(x) \sum_{i=1}^{\infty} \E\left[\widetilde{\Theta}_i^{\beta}(t)\,{\bf 1}_{\{\theta >0\}} \right] ,
\eeam 
uniformly for $t \in \Lambda_T$, for any fixed  $T \in \Lambda$. From \eqref{eq.KPY.3.9} and \eqref{eq.KPY.3.10}, we see immediately the effect of the dependence between insurance and financial risks through the distributions $Q_i$, with $i \in \bbn$. Due to Assumption \ref{ass.KPY.3.1}, recall Remark \ref{rem.KPY.3.1}, it holds either that $\bV(x) = o[\overline{Q}_i(x)]$,  or $\bV(x) \asymp \overline{Q}_i(x)$, for $i \in \bbn$. In the first case, the first term of the right member of \eqref{eq.KPY.3.9} (or, of \eqref{eq.KPY.3.10}, respectively), means that the financial risk has asymptotically negligible effect on the decay rate of the probability from \eqref{eq.KPY.1.7}. \cite{guo:2022}, \cite{xu:peng:zou:2025}, in an one-dimensional renewal risk model, considered only the case $\bV(x) \asymp \overline{Q}_i(x)$, that excludes some interesting cases of strongly dependent insurance and financial risks, through ${\bf X}^{(i)}$, $e^{-Y_i}$. 
\ere

\section{Argumentation of Theorem \ref{th.KPY.3.1}} \label{sec.KPT29U.4}

In this Section we present the proof of Theorem \ref{th.KPY.3.1}, and of Corollary, after some preliminary lemmas. We note that the $\Theta_i(t)$, $\widetilde{\Theta}_i(t)$, with $i \in \bbn$, from relation \eqref{eq.KPY.3.6} take the form
\beam \label{eq.KPY.4.1}
&&\Theta_i(t)\,{\bf 1}_{\{\theta=0\}} = \exp\left\{-B(\tau_i + \theta) - \sum_{j=1}^{i-1} Y_j\right\}\,{\bf 1}{\{\tau_j \leq t\}}\,{\bf 1}_{\{\theta=0\}}\,, \\[2mm] \notag
&& \widetilde{\Theta}_i(t)\,{\bf 1}_{\{\theta>0\}} = \exp\left\{-B(\tau_i + \theta)+Y_i - \sum_{j=1}^{N(\tau_i+\theta)-} Y_j\right\}\,{\bf 1}{\{\tau_j +\theta \leq t\}}
\,{\bf 1}_{\{\theta> 0\}}\,.
\eeam
Hence, for all $i \in \bbn$, the ${\bf X}^{(i)}$ and $\Theta_i(t)\,{\bf 1}_{\{\theta=0\}}$ are independent, and furthermore the $\widehat{\bf X}^{(i)}$
and $\widetilde{\Theta}_i(t)\,{\bf 1}_{\{\theta>0\}}$ are independent, by Assumption \ref{ass.KPY.1.2}. During the proofs of this section, as also of Subsection \ref{subsec.KPT29U.5.2}, when we say $t \in \Lambda_T$, we mean $t \in \Lambda_T$ for any fixed $T \in \Lambda$.

The first lemma provides the presence of multivariate linear single big jump principle on finite number randomly weighted sums. In this lemma we use generalized moment conditions from Assumption \ref{ass.KPY.3.1}, while the restriction to class $\mathcal{P_D}_A$ is not needed.

\ble \label{lem.KPY.4.1}
Let $A \in \mathscr{R}$ some fixed set. We suppose that Assumptions \ref{ass.KPY.1.1}, \ref{ass.KPY.1.2} hold, the $\{ {\bf X}^{(i)}\,,\;i \in \bbn \}$ are $QAI_A$, with $F \in \mathcal{C}_A$, and 
$G_i \in \mathcal{C}_A$, with $i \in \bbn$, with $G_i(x\,A) \asymp G_j(x\,A)$, $i,\,j \in \bbn$, $i \neq j$. Let, for any $n \in \bbn$, there exists some 
\beao
p > J_{F_A}^+ \vee  J_{G_A^{(1)}}^+ \,,
\eeao 
such that it holds $M_p^{(i)} < \infty$, for $i = 1,\,\ldots,\,n$.  Then, it holds
\beam \label{eq.KPY.4.2}
&&\PP \left( \sum_{i=1}^{n} {\bf X}^{(i)}\,\Theta_i(t) \in x\,A \right) 
\sim \sum_{i=1}^{n} \PP \left( {\bf X}^{(i)}\,\Theta_i(t)\,{\bf 1}_{\{\theta=0\}} \in x\,A \right)  \\[2mm] \notag
&&+  \sum_{i=1}^{n} \PP \left( \widehat{\bf X}^{(i)}\,\widetilde{\Theta}_i(t)\,{\bf 1}_{\{\theta>0\}}\in x\,A \right) = \sum_{i=1}^{n} \PP \left( {\bf X}^{(i)}\,\Theta_i(t) \in x\,A \right)\,,
\eeam
uniformly for $t \in  \Lambda_T$, for any fixed $T \in \Lambda$.
\ele

\pr~
At first we shall show that the sequence $\left\{ \widehat{\bf X}^{(i)} \,,\; i \in\bbn \right\}$ is $QAI_A$. It is enough to show that the sequence $\left\{\widehat{ X}_A^{(i)}\,,\; i \in\bbn \right\}$ (recall relation \eqref{eq.KPY.3.a}) is $QAI$. Since $F_A \in \mathcal{C}$ and $M_p^{(i)} < \infty$ for any $i \in \bbn $, in combination with the fact that the sequence $\left\{X_A^{(i)}\,,\; i \in\bbn \right\}$ is $QAI$, by \cite[Th. 2.2]{li:2013} we obtain that the sequence $\left\{ \widehat{ X}_A^{(i)}=X_A^{(i)}\,e^{-Y_i}\,,\; i \in\bbn \right\}$ is $QAI$. 

Next, it holds
\beam \label{eq.KPY.4.4}
&&\PP \left( \sum_{i=1}^{n} {\bf X}^{(i)}\,\Theta_i(t) \in x\,A \right) 
  \\[2mm] \notag
&&=\PP \left( \sum_{i=1}^{n} {\bf X}^{(i)}\,\Theta_i(t)\,{\bf 1}_{\{\theta=0\}}  \in x\,A \right)+\PP \left( \sum_{i=1}^{n}  \widehat{\bf X}^{(i)}\,\widetilde{\Theta}_i(t)\,{\bf 1}_{\{\theta>0\}} \in x\,A \right) \,,
\eeam
for all $t \in  \Lambda_T$. We have to show that for some $p > J_{F_A}^+ \vee J_{G_A^{(1)}}^+ $ it holds
\beam \label{eq.KPY.4.5}
\E \left[ \Theta_i^p(t) \,{\bf 1}_{\{\theta=0\}} \right] < \infty\,, \quad  \E \left[ \widetilde{\Theta}_i^p(t) \,{\bf 1}_{\{\theta>0\}} \right] < \infty\,, 
\eeam
for all $t \in  \Lambda_T$.  

Indeed, for any $t \in  \Lambda_T$ we obtain
\beao
&&\E \left[ \left(\Theta_i(t) \right)^p\,{\bf 1}_{\{\theta=0\}} \right] = \E \left[ \left(\exp\left\{-B(\tau_i ) - \sum_{j=1}^{i-1} Y_j\right\}\,{\bf 1}{\{\tau_j \leq t\}}\right)^p\,{\bf 1}_{\{\theta=0\}} \right] \\[2mm]
&& =\PP(\theta=0)\,\int_0^t \E \left[ \left( e^{-B(s ) - \sum_{j=1}^{i-1} Y_j}\,{\bf 1}{\{\tau_j \leq t\}}\right)^p \right]\,\PP(\tau_i \in ds) \\[2mm]
&& \leq \E\left[e^{ -p\, \sum_{j=1}^{i-1} Y_j}\right] \,\int_0^t \E \left[ e^{-p\,B(s ) }\right]\,\PP(\tau_i \in ds) \\[2mm]
&& \leq \E\left[e^{ -p\, (i-1)\,\widehat{Y}}\right] \,\int_0^t e^{s\,\phi_B(p)}\,\PP(\tau_i \in ds) \leq \left\{\E\left[\left(e^{ -\widehat{Y}}\right)^p\right]\right\}^{i-1} \,\E\left[e^{\tau_i\,\phi_B(p)}\,{\bf 1}_{\{\tau_i \leq T\}} \right] < \infty\,, 
\eeao 
where we denote 
\beao
\widehat{Y}:= \bigwedge_{j=1}^{i-1} Y_j\,,
\eeao 
and at the fourth step we used \eqref{eq.KPY.1.3} and the last step follows from the inequalities $-\infty < \phi_B(p) < \infty$, hence the last expectation is bounded either by $e^{T\,\phi_B(p)}$, when $\phi_B(p) > 0$, or by unity, when $\phi_B(p) \leq 0$. The last relation provides the first inequality of \eqref{eq.KPY.4.5}.

The second inequality of \eqref{eq.KPY.4.5} follows by quite similar way, through integration with respect to $\PP(\tau_i + \theta \in ds)$.

From the first relation of \eqref{eq.KPY.4.5}, and since $F \in \mathcal{C}_A$, we find, for  any $t \in  \Lambda_T$, the following statements:
\begin{enumerate}
\item
The $X_A^{(i)}\,\Theta_i(t)\,{\bf 1}_{\{\theta=0\}} $ follows a distribution from class $\mathcal{C}$, see \cite[Th.3.4 (ii)]{cline:samorodnitsky:1994} or \cite[Prop. 5.3 (iv)]{leipus:siaulys:konstantinides:2023}. Recall that the $X_A^{(i)}$ and $\Theta_i(t)\,{\bf 1}_{\{\theta=0\}} $ are independent. Hence the distribution of the product  $ {\bf X}^{(i)}\,\Theta_i(t) $ belongs to class $\mathcal{C}_A$
\item
The sequence $\left\{ {\bf X}^{(i)}\,\Theta_i(t)\,{\bf 1}_{\{\theta=0\}}\,,\; i \in\bbn \right\}$ is $QAI_A$, due to \cite[Th. 2.2]{li:2013}, that is valid because of the first relation of \eqref{eq.KPY.4.5}, in combination with $F_A \in \mathcal{C}$.  
\end{enumerate}

Therefore, by \cite[Th. 4.1 (i)]{konstantinides:passalidis:2024g}, we obtain
\beam \label{eq.KPY.4.6}
\PP \left( \sum_{i=1}^{n} {\bf X}^{(i)}\,\Theta_i(t)\,{\bf 1}_{\{\theta=0\}}  \in x\,A \right) \sim \sum_{i=1}^{n} \PP \left( {\bf X}^{(i)}\,\Theta_i(t)\,{\bf 1}_{\{\theta=0\}}  \in x\,A \right)\,, 
\eeam
for all $t \in  \Lambda_T$.

In a similar way, using the second inequality of \eqref{eq.KPY.4.5}, for all $t \in  \Lambda_T$, the sequence 
\beao
\left\{ \widehat{\bf X}^{(i)}\,\widetilde{\Theta}_i(t)\,{\bf 1}_{\{\theta>0\}}\,,\; i \in\bbn \right\}\,,
\eeao 
is $QAI_A$, and their distributions belong to class $\mathcal{C}_A$. Applying \cite[Th. 4.1 (i)]{konstantinides:passalidis:2024g} we find that it holds
\beam \label{eq.KPY.4.7}
\PP \left( \sum_{i=1}^{n} \widehat{\bf X}^{(i)}\,\widetilde{\Theta}_i(t)\,{\bf 1}_{\{\theta>0\}}  \in x\,A \right) \sim \sum_{i=1}^{n} \PP \left( \widehat{\bf X}^{(i)}\,\widetilde{\Theta}_i(t)\,{\bf 1}_{\{\theta>0\}} \in x\,A \right)\,, 
\eeam
for all $t \in  \Lambda_T$.

From \eqref{eq.KPY.4.4}, \eqref{eq.KPY.4.6} and \eqref{eq.KPY.4.7}, we conclude that \eqref{eq.KPY.4.2} holds uniformly for $t \in  \Lambda_T$.     
~\halmos

\ble \label{lem.KPY.4.2}
Let $A \in \mathscr{R}$ some fixed set. We suppose that Assumptions \ref{ass.KPY.1.1}, \ref{ass.KPY.1.2}, \ref{ass.KPY.3.1}, \ref{ass.KPY.3.2} are valid.  Then, it holds
\beam \label{eq.KPY.4.8}
\lim_{N \to \infty} \limsup \sup_{t\in \Lambda_T} \dfrac{\PP \left( \sum_{i=N+1}^{\infty} {\bf X}^{(i)}\,\Theta_i(t) \in x\,A \right)}{\PP \left( {\bf X}^{(1)}\,\Theta_1(t)\,{\bf 1}_{\{\theta=0\}}  \in x\,A \right) +\PP \left( \widehat{\bf X}^{(1)}\,\widetilde{\Theta}_1(t)\,{\bf 1}_{\{\theta>0\}} \in x\,A \right)} = 0\,,
\eeam
and 
\beam \label{eq.KPY.4.9}
\lim_{N \to \infty} \limsup \sup_{t\in \Lambda_T} \dfrac{\sum_{i=N+1}^{\infty} \PP \left( {\bf X}^{(i)}\,\Theta_i(t) \in x\,A \right)}{\PP \left( {\bf X}^{(1)}\,\Theta_1(t)\,{\bf 1}_{\{\theta=0\}}  \in x\,A \right) +\PP \left( \widehat{\bf X}^{(1)}\,\widetilde{\Theta}_1(t)\,{\bf 1}_{\{\theta>0\}} \in x\,A \right)} = 0\,,
\eeam
for any fixed $T \in  \Lambda$.
\ele

\pr~
At first, taking into account that $F \in \mathcal{C}_A$, by \eqref{eq.KPY.4.5}, recalling that Assumption \ref{ass.KPY.3.2} is stronger than the moment conditions that we assumed in Lemma \ref{lem.KPY.4.1}, from \cite[Th. 3.3 (iv)]{cline:samorodnitsky:1994} we obtain that it holds
\beam \label{eq.KPY.4.10}
&&\PP \left( {\bf X}^{(i)}\,\Theta_i(t)\,{\bf 1}_{\{\theta=0\}}  \in x\,A \right)= \PP \left( X_A^{(i)}\,\Theta_i(t)\,{\bf 1}_{\{\theta=0\}} > x \right) \\[2mm] \notag
&& \asymp \PP \left( X_A > x \right) = \PP \left( {\bf X}^{(1)} \in x\,A \right)\,,
\eeam
for any $i \in \bbn$ and all $t \in  \Lambda_T$, and similarly, since $G_i \in \mathcal{C}_A$ and $G_i(x\,A) \asymp G_1(x\,A)$, for any $i \in \bbn$ it holds
 \beam \label{eq.KPY.4.11}
\PP \left( \widehat{\bf X}^{(i)}\,\widetilde{\Theta}_i(t)\,{\bf 1}_{\{\theta>0\}} \in x\,A \right)\asymp  \PP \left( \widehat{\bf X}^{(1)} \in x\,A \right)\,,
\eeam
for all $t \in  \Lambda_T$. From \eqref{eq.KPY.4.10} is implied that the 
\beao
X_A^{(i)}\,\Theta_i(t)\,{\bf 1}_{\{\theta=0\}}\,,
\eeao 
with $i \in \bbn$ and all $t \in  \Lambda_T$, have Matuszewska index equal to $J_{F_A}^+$, while from \eqref{eq.KPY.4.11} is implied that the 
\beao
\widehat{ X}_A^{(i)}\,\widetilde{\Theta}_i(t)\,{\bf 1}_{\{\theta>0\}}\,,
\eeao 
with $i \in \bbn$ and all $t \in  \Lambda_T$, have Matuszewska index equal to $J_{G_A^{(1)}}^+$.

Further, it holds
\beam \label{eq.KPY.4.12}
&&\PP \left( \sum_{i=N+1}^{\infty} {\bf X}^{(i)}\,\Theta_i(t) \in x\,A \right) 
  \\[2mm] \notag
&&=\PP \left( \sum_{i=N+1}^{\infty} {\bf X}^{(i)}\,\Theta_i(t)\,{\bf 1}_{\{\theta=0\}}  \in x\,A \right)+\PP \left( \sum_{i=N+1}^{n}  \widehat{\bf X}^{(i)}\,\widetilde{\Theta}_i(t)\,{\bf 1}_{\{\theta>0\}} \in x\,A \right) \,,
\eeam 
for all $t \in  \Lambda_T$. For the fixed $\gamma > 1$ from Assumption \ref{ass.KPY.3.2}, we can find some large enough $n_0 \bbn$, such that for any $N \geq n_0$ it holds $\sum_{i=N+1}^{\infty} i^{-\gamma} \leq 1$. Then we obtain
\beam \label{eq.KPY.4.13} \notag
&&\PP \left( \sum_{i=N+1}^{\infty} {\bf X}^{(i)}\,\Theta_i(t)\,{\bf 1}_{\{\theta=0\}}  \in x\,A \right)  =\PP \left( \sup_{{\bf p}\in I_A} {\bf p}^{\top} \left[\sum_{i=N+1}^{\infty} {\bf X}^{(i)}\,\Theta_i(t)\,{\bf 1}_{\{\theta=0\}} \right] > x \right) \\[2mm] \notag
&&\leq \PP \left( \sum_{i=N+1}^{\infty} X_A^{(i)}\,\Theta_i(t)\,{\bf 1}_{\{\theta=0\}} > x \right) \leq \PP \left( \sum_{i=N+1}^{\infty} X_A^{(i)}\,\Theta_i(t)\,{\bf 1}_{\{\theta=0\}} > x\,\sum_{i=N+1}^{\infty} i^{-\gamma} \right) \\[2mm] \notag
&&\leq \PP \left( \bigcup_{i=N+1}^{\infty} \left\{X_A^{(i)}\,\Theta_i(t)\,{\bf 1}_{\{\theta=0\}} > x\,i^{-\gamma} \right\} \right) \leq \sum_{i=N+1}^{\infty} \PP \left( X_A^{(i)}\,\Theta_i(t)\,{\bf 1}_{\{\theta=0\}} > x\,i^{-\gamma} \right) \,,\\
\eeam 
for all $t \in  \Lambda_T$, where at the second step we used the fact that the sequences of events 
\beao
\left\{ \left( \sum_{i=N+1}^{k} {\bf X}^{(i)}\,\Theta_i(t)\,{\bf 1}_{\{\theta=0\}}  \in x\,A \right)\; k \in \bbn \right\}\,,\; \left\{ \left( \sum_{i=N+1}^{k} X_A^{(i)}\,\Theta_i(t)\,{\bf 1}_{\{\theta=0\}} > x \right)\,,\; k \in \bbn \right\}
\eeao 
are increasing. 

We define
\beam \label{eq.KPY.ex.1}
R^{+}(T) : = \exp\left\{ - \inf_{0\leq s \leq T} B(s) - \inf_{0\leq s \leq T} J(s) \right\}\,.
\eeam
Then, the following inequality holds almost surely for all $t \in \Lambda_T$
\beam \label{eq.KPY.ex.2}
\Theta_i(t)\,{\bf 1}_{\{\theta =0\}} \leq  R^{+}(T)\,{\bf 1}_{\{\theta =0\}}\,{\bf 1}_{\{\tau_i \leq t \}}\,.
\eeam
For each term of the last sum in \eqref{eq.KPY.4.13}, from \eqref{eq.KPY.ex.2}, for all $t \in \Lambda_T$ it holds
\beam \label{eq.KPY.ex.3}
&&\PP(X_A^{(i)}\,\Theta_i(t)\,{\bf 1}_{\{\theta =0\}} > x\,i^{-\gamma} ) \leq  \PP(X_A^{(i)}\,R^{+}(T)\,{\bf 1}_{\{\theta =0\}} > x\,i^{-\gamma}\,,\;\tau_i \leq t )  \\[2mm] \notag
&& \leq  \PP(X_A^{(i)}\,R^{+}(T) > x\,i^{-\gamma} )\,\PP(\theta=0)\,\PP(\tau_i \leq T) \,.
\eeam 
Following the argumentation in the proof of \cite[Lem. 4.4]{guo:2022}, and applying the generalized Potter inequalities for the $q^*$ in Assumption \ref{ass.KPY.3.2}, see \cite[Prop. 2.2.1]{bingham:goldie:teugels:1987}, we find that for any large enough $x>0$, and any $i \in \bbn$ it holds
\beam \label{eq.KPY.ex.4}
\PP(X_A^{(i)}\,R^{+}(T)> x\,i^{-\gamma} ) \leq C\,i^{\gamma\,q^*} \PP(X_A^{(i)}> x ) \,,
\eeam
for some constant $C = C\left( \E\left[ \left(R^{+}(T) \right)^{\gamma\,q^*} \right]\right)$, 
where 
\beao
\E\left[ \left(R^{+}(T)\right)^{\gamma\,q^*} \right]=\E\left[  \exp\left\{ - \gamma\,q^*\,\left[\inf_{0\leq s \leq T} B(s) - \inf_{0\leq s \leq T} J(s)\right] \right\} \right] < \infty\,,
\eeao
that is implied by Assumption \ref{ass.KPY.3.1}. Hence, by \eqref{eq.KPY.4.13}, \eqref{eq.KPY.ex.3} and \eqref{eq.KPY.ex.4}, for all $t \in \Lambda_T$, $N \geq n_0$, we obtain
\beam \label{eq.KPY.ex.5}
&&\PP\left(\sum_{i=N+1}^{\infty} {\bf X}^{(i)}\,\Theta_i(t)\,{\bf 1}_{\{\theta =0\}} \in x\,A \right) \lesssim C\,\PP( {\bf X} \in x\,A) \sum_{i=N+1}^{\infty} \,i^{\gamma\,q^*}\PP(\tau_i \leq T) \\[2mm] \notag
&& = \dfrac C{\gamma\,q^* +1}\,\PP( {\bf X} \in x\,A) \E\left[ N^{\gamma\,q^* +1}(T)\,{\bf 1}_{\{N(T) > N+1 \}} \right] \,.
\eeam

Hence, from \eqref{eq.KPY.ex.5} we obtain that for any $\vep >0$, there exists some $n^* \in \bbn$ such that for any $N > n^*\vee n_0$, it holds
\beam \label{eq.KPY.4.16}
&&\PP \left( \sum_{i=N+1}^{\infty} {\bf X}^{(i)}\,\Theta_i(t)\,{\bf 1}_{\{\theta=0\}}  \in x\,A \right)  \\[2mm] \notag
&&\lesssim \vep\,\PP({\bf X}^{(1)} \in x\,A)  \lesssim \vep\,C''\,\PP({\bf X}^{(1)}\,\Theta_1(t)\,{\bf 1}_{\{\theta=0\}} \in x\,A) \,,
\eeam
for all $t \in \Lambda_T$, where the $C'' \in (0,\,\infty)$ follows by \eqref{eq.KPY.4.10}.

For the second probability on the right member of \eqref{eq.KPY.4.12},  following similar steps, as these for \eqref{eq.KPY.4.13}, \eqref{eq.KPY.ex.5}  we find 
\beam \label{eq.KPY.4.17}
&&\PP \left( \sum_{i=N+1}^{\infty} \widehat{\bf X}^{(i)}\,\widetilde{\Theta}_i(t)\,{\bf 1}_{\{\theta>0\}} \in x\,A \right) \\[2mm] \notag
&& \lesssim \dfrac C{\gamma\,q^* +1}\,\PP(\widehat{\bf X}^{(1)}\in x\,A)\,\E\left[N^{\gamma\,q^* +1}(T)\,{\bf 1}_{\{N(T) > N+1\}} \right] \,,
\eeam
for all $t \in  \Lambda_T$, where now the constant $C$, is derived for the term $\PP(X_A^{(i)}\,e^{-\widehat{Y}} > x ) \geq 
\PP(X_A^{(i)}\,e^{-Y_i} > x ) $, for any $x > 0$. 

Hence, from \eqref{eq.KPY.4.17}, for any $\vep>0$, there exists some $n' \in\bbn$, such that for all $N >n_0\vee n'$ it holds
\beam \label{eq.KPY.4.19}
&&\PP \left( \sum_{i=N+1}^{\infty} \widehat{\bf X}^{(i)}\,\widetilde{\Theta}_i(t)\,{\bf 1}_{\{\theta>0\}} \in x\,A \right) \\[2mm] \notag
&& \lesssim \vep\,\PP(\widehat{\bf X}^{(1)}\in x\,A) \lesssim \vep\,C''\,\PP(\widehat{\bf X}^{(1)}\,\widetilde{\Theta}_1(t)\,{\bf 1}_{\{\theta>0\}} \in x\,A)\,,
\eeam
for all $t \in  \Lambda_T$, where at the last step on $C'' \in (0,\,\infty)$ we used \eqref{eq.KPY.4.11}.

From \eqref{eq.KPY.4.12}, \eqref{eq.KPY.4.16}, \eqref{eq.KPY.4.19}, and the arbitrary choice of $\vep>0$, we obtain \eqref{eq.KPY.4.8}.

Relation  \eqref{eq.KPY.4.9} can be proved similarly and easier, since we do not need the sum $\sum_{i=N+1}^{\infty} i^{-\gamma} $, as in \eqref{eq.KPY.4.13}.
~\halmos

The following lemma plays crucial role for the argumentation of Theorem \ref{th.KPY.3.1} and extends Lemma \ref{lem.KPY.4.1}, uniformly for $n \in \bbn$, under some additional conditions.

\ble \label{lem.KPY.4.3}
Let $A \in \mathscr{R}$ be some fixed set. Under Assumptions  \ref{ass.KPY.1.1}, \ref{ass.KPY.1.2}, \ref{ass.KPY.3.1}, \ref{ass.KPY.3.2}, then for any fixed $T \in \Lambda$, it holds  
\beam \label{eq.KPY.4.20} \notag
&&\lim \sup_{n\in \bbn} \left|\dfrac{\PP \left( \sum_{i=1}^{n} {\bf X}^{(i)}\,\Theta_i(t) \in x\,A \right)}{ \sum_{i=1}^{n} \PP \left({\bf X}^{(i)}\,\Theta_i(t)\,{\bf 1}_{\{\theta=0\}}  \in x\,A \right)+  \sum_{i=1}^{n}
\PP \left( \widehat{\bf X}^{(1)}\,\widetilde{\Theta}_1(t)\,{\bf 1}_{\{\theta>0\}} \in x\,A \right)} -1\right|\\[2mm] 
&& =\lim \sup_{n\in \bbn} \left|\dfrac{\PP \left( \sum_{i=1}^{n} {\bf X}^{(i)}\,\Theta_i(t) \in x\,A \right)}{ \sum_{i=1}^{n} \PP \left({\bf X}^{(i)}\,\Theta_i(t) \in x\,A \right)} -1\right|=0\,,
\eeam
uniformly for $t \in  \Lambda_T$.
\ele

\pr~ 
We show only the first relation of \eqref{eq.KPY.4.20}, since the second follows immediately from it. Let consider a fixed, large $N \in \bbn$. Then by Lemma \ref{lem.KPY.4.1} we obtain
\beam \label{eq.KPY.4.21} \notag
\lim \sup_{N\geq n \in \bbn} \left|\dfrac{\PP \left( \sum_{i=1}^{n} {\bf X}^{(i)}\,\Theta_i(t) \in x\,A \right)}{ \sum_{i=1}^{n} \PP \left({\bf X}^{(i)} \Theta_i(t) {\bf 1}_{\{\theta=0\}}  \in x A \right)+  \sum_{i=1}^{n}
\PP \left( \widehat{\bf X}^{(1)} \widetilde{\Theta}_1(t) {\bf 1}_{\{\theta>0\}} \in x A \right)} -1\right| =0 ,
\eeam 
uniformly for $t \in  \Lambda_T$. It remains to show the first relation of \eqref{eq.KPY.4.20} for all $n> N$.

For the lower bound for all $n> N$, we find
\beao
&&\PP \left( \sum_{i=1}^{n} {\bf X}^{(i)}\,\Theta_i(t) \in x\,A \right) \geq \PP \left( \sum_{i=1}^{N} {\bf X}^{(i)}\,\Theta_i(t) \in x\,A \right) \\[2mm]
&&\sim \sum_{i=1}^{N} \PP \left( {\bf X}^{(i)}\,\Theta_i(t)\,{\bf 1}_{\{\theta=0\}} \in x\,A \right) + \sum_{i=1}^{N} \PP \left( \widehat{\bf X}^{(i)}\,\widetilde{\Theta}_i(t)\,{\bf 1}_{\{\theta>0\}} \in x\,A \right) \\[2mm]
&&\geq \left(\sum_{i=1}^{n} - \sum_{i=N+1}^{\infty}\right) \PP \left( {\bf X}^{(i)}\,\Theta_i(t)\,{\bf 1}_{\{\theta=0\}} \in x\,A \right) \\[2mm] 
&&+\left(\sum_{i=1}^n - \sum_{i=N+1}^{\infty} \right) \PP \left( \widehat{\bf X}^{(i)}\,\widetilde{\Theta}_i(t)\,{\bf 1}_{\{\theta>0\}} \in x\,A \right)\,,
\eeao
uniformly for $t \in  \Lambda_T$, where at the second step we applied Lemma \ref{lem.KPY.4.1}. From Lemma \ref{lem.KPY.4.2}, recall relations \eqref{eq.KPY.4.16} and \eqref{eq.KPY.4.19}, for any $\vep > 0$, we can find a large enough $N>0$, such that for all $n> N$, it holds
\beam \label{eq.KPY.4.22} \notag
&&\PP \left( \sum_{i=1}^{n} {\bf X}^{(i)}\,\Theta_i(t) \in x\,A \right) \gtrsim  \sum_{i=1}^{n} \PP \left( {\bf X}^{(i)}\,\Theta_i(t)\,{\bf 1}_{\{\theta=0\}} \in x\,A \right) \\[2mm]
&&-\vep\, \PP \left( {\bf X}^{(1)}\,\Theta_1(t)\,{\bf 1}_{\{\theta=0\}} \in x\,A \right) + \sum_{i=1}^{n} \PP \left( \widehat{\bf X}^{(i)}\,\widetilde{\Theta}_i(t)\,{\bf 1}_{\{\theta>0\}} \in x\,A \right) \\[2mm] \notag
&&-\vep\, \PP \left( \widehat{\bf X}^{(1)}\,\widetilde{\Theta}_1(t)\,{\bf 1}_{\{\theta>0\}} \in x\,A \right)  \\[2mm] \notag
&&\geq (1-\vep) \Bigg(\sum_{i=1}^{n} \PP \left( {\bf X}^{(i)}\,\Theta_i(t)\,{\bf 1}_{\{\theta=0\}} \in x\,A \right) +\sum_{i=1}^n \PP \left( \widehat{\bf X}^{(i)}\,\widetilde{\Theta}_i(t)\,{\bf 1}_{\{\theta>0\}} \in x\,A \right) \Bigg)\,,
\eeam
uniformly for $t \in  \Lambda_T$.

From the other hand side, for any $\vep >0$, for all $n > N$, we obtain
\beam \label{eq.KPY.4.23} \notag
&&\PP \left( \sum_{i=1}^{n} {\bf X}^{(i)}\,\Theta_i(t) \in x\,A \right) \leq
\PP \left( \sum_{i=1}^{n}  X_{A}^{(i)}\,\Theta_i(t) >x \right) \leq \PP \left( \sum_{i=1}^{\infty}  X_{A}^{(i)}\,\Theta_i(t) >x \right) \\[2mm] 
&&= \PP \left( \sum_{i=1}^{\infty}  X_{A}^{(i)}\,\Theta_i(t)\,{\bf 1}_{\{\theta=0\}} >x \right) +\PP \left( \sum_{i=1}^{\infty} \widehat{X}_A^{(i)}\,\widetilde{\Theta}_i(t)\,{\bf 1}_{\{\theta>0\}} >x \right) \\[2mm]\notag
&&\sim \sum_{i=1}^{\infty} \PP \left(  X_{A}^{(i)}\,\Theta_i(t)\,{\bf 1}_{\{\theta=0\}} >x \right) + \sum_{i=1}^{\infty} \PP \left(\widehat{X}_A^{(i)}\,\widetilde{\Theta}_i(t)\,{\bf 1}_{\{\theta>0\}} >x \right)\leq \left(\sum_{i=1}^n + \sum_{i=N+1}^{\infty}\right) \\[2mm]  \notag
&& \PP \left( {\bf X}^{(i)}\,\Theta_i(t)\,{\bf 1}_{\{\theta=0\}} \in x\,A \right) + \left(\sum_{i=1}^n + \sum_{i=N+1}^{\infty} \right) \PP \left(\widehat{\bf X}^{(i)}\,\widetilde{\Theta}_i(t)\,{\bf 1}_{\{\theta>0\}} \in x\,A \right) \\[2mm] \notag
&&\leq(1+\vep)\,\left[\sum_{i=1}^n \PP \left({\bf X}^{(i)}\,\Theta_i(t)\,{\bf 1}_{\{\theta=0\}} \in x\,A \right)  + \sum_{i=1}^n  \PP \left(\widehat{\bf X}^{(i)}\,\widetilde{\Theta}_i(t)\,{\bf 1}_{\{\theta>0\}} \in x\,A \right)\right] 
\,,
\eeam
uniformly for $t \in  \Lambda_T$, where at the first step we applied \cite[Prop. 2.4]{konstantinides:passalidis:2024g}, and at the fourth step we used \cite[Th. 2]{yi:chen:su:2011}, which we explain its application at the end of the proof. At the last step we applied Lemma \ref{lem.KPY.4.2}. Hence, by \eqref{eq.KPY.4.22},  \eqref{eq.KPY.4.23}, due to the arbitrary choice of $\vep>0$, we find
\beam \label{eq.KPY.4.24} \notag
\limsup_{N < n} \left|\dfrac{\PP \left( \sum_{i=1}^{n} {\bf X}^{(i)}\,\Theta_i(t) \in x\,A \right)}{\sum_{i=1}^n \PP \left({\bf X}^{(i)}\,\Theta_i(t)\,{\bf 1}_{\{\theta=0\}} \in x\,A \right)  + \sum_{i=1}^n  \PP \left(\widehat{\bf X}^{(i)}\,\widetilde{\Theta}_i(t)\,{\bf 1}_{\{\theta>0\}} \in x\,A \right)} - 1 \right| =0 \,,\\
\eeam
uniformly for $t \in  \Lambda_T$. From relations \eqref{eq.KPY.4.23} and  \eqref{eq.KPY.4.24} we obtain that the first limit of  \eqref{eq.KPY.4.20} tends to zero. 

We close the proof with the justification of the application of \cite[Th. 2]{yi:chen:su:2011}. If $J_{F_A}^+ <1$, then \cite[Th. 2]{yi:chen:su:2011} in order to be applied we must show that 
for $q^*$ from Assumption \ref{ass.KPY.3.2}, and for all $t \in  \Lambda_T$, it holds
\beam \label{eq.KPY.4.14}
\sum_{i=1}^{\infty} i^{\gamma\,q^*}\,\E\left[\Theta_i^{p_1}(t)\,{\bf 1}_{\{\theta=0\}} \bigvee \Theta_i^{p_2}(t)\,{\bf 1}_{\{\theta=0\}}\right] < \infty \,.
\eeam
Indeed, it is enough to show only  
\beam \label{eq.KPY.4.15}
\sum_{i=1}^{\infty} i^{\gamma\,q^*} \,\E\left[\Theta_i^{p_2}(t)\,{\bf 1}_{\{\theta=0\}} \right] < \infty \,,
\eeam
since the corresponding convergence with index $p_1$ follows similarly.

From the independence between $\{N(t)\,,\; t\geq 0\}$ and $\theta$, and the independence from all the other sources of randomness, recall Assumption \ref{ass.KPY.1.2}, we obtain 
\beao
&&\sum_{i=1}^{\infty} i^{\gamma\,q^*}\,\E\left[ \Theta_i^{p_2}(t)\,{\bf 1}_{\{\theta=0\}}\right]  =\PP(\theta=0)\,\sum_{i=1}^{\infty} i^{\gamma\,q^*} \,\E\left[\exp\left\{-B(\tau_i)-\sum_{j=1}^{i-1}Y_j\right\}\,{\bf 1}_{\{\tau_i \leq t\}} \right]^{p_2} \\[2mm] \notag
&& \leq \sum_{i=1}^{\infty} i^{\gamma\,q^*} \,\int_0^t \E\left[\left(e^{-\sum_{j=1}^{i-1}Y_j}\right)^{p_2}\right]\,\E\left[\left(e^{-B(s)} \right)^{p_2}\right]\,\PP(\tau_i \in ds)  \\[2mm] \notag
&& \leq \sum_{i=1}^{\infty} i^{\gamma\,q^*} \,\left(\E\left[\left(e^{-\widehat{Y}}\right)^{p_2}\right]\right)^{i-1}\,\int_0^t e^{s\,\phi(p_2)} \,\PP(\tau_i \in ds) \\[2mm] \notag
&& \leq C^* \sum_{i=1}^{\infty} i^{\gamma\,q^*} \,\left(\E\left[\left(e^{-\widehat{Y}}\right)^{p_2}\right]\right)^{i-1}\,\PP(\tau_i \leq t) \leq \dfrac{C^*}{\gamma\,q^* +1}\,\E[(N(t))^{\gamma\,q^*+1}] < \infty\,, 
\eeao
for all $t \in  \Lambda_T$, where $C^*= e^{T\,\phi_B(p_2)}$, when $\phi_B(p_2)>0$ and $C^*=1$, when $\phi_B(p_2)\leq 0$. The last step follows by Assumption \ref{ass.KPY.3.2}. If $J_{G_A^{(1)}}^+ <1$, the same theorem can be applied on the second sum through 
\beam \label{eq.KPY.4.18}
\sum_{i=1}^{\infty} i^{\gamma\,q^*}\,\E\left[\widetilde{\Theta}_i^{q_1}(t)\,{\bf 1}_{\{\theta>0\}} \bigvee \widetilde{\Theta}_i^{q_2}(t)\,{\bf 1}_{\{\theta>0\}}\right] < \infty\,,
\eeam
for all $t \in  \Lambda_T$. It is enough to show only
\beao
\sum_{i=1}^{\infty} i^{\gamma\,q^*}\,\E\left[\widetilde{\Theta}_i^{q_2}(t)\,{\bf 1}_{\{\theta>0\}}\right] < \infty\,,
\eeao
for all $t \in  \Lambda_T$, since the corresponding convergence with $q_1$, can be proved through similar arguments.

Following the line of proof for \eqref{eq.KPY.4.15}, taking into consideration that since $\theta>0$, we obtain $N[(\tau_i+ \theta)-] \geq i$, and $N_{\Theta}(t) \leq N(t)$ almost surely, for all $t \in  \Lambda_T$, we find
\beao
&&\sum_{i=1}^{\infty} i^{\gamma\,q^*}\,\E\left[\widetilde{\Theta}_i^{q_2}(t)\,{\bf 1}_{\{\theta>0\}} \right] =\PP(\theta >0)\,\sum_{i=1}^{\infty} i^{\gamma\,q^*}\,\E\left[\left(e^{-B(\tau_i + \theta) +Y_i - \sum_{j=1}^{N[(\tau_i+\theta)-]}}\, {\bf 1}_{\{\tau_i+\theta \leq t\}}\right)^{q_2} \right] \\[2mm] 
&&\leq \sum_{i=1}^{\infty} i^{\gamma\,q^*}\,\int_0^t\,\E\left[\left(e^{-B(s)}\right)^{p_2} \right]\,\E\left[\left(e^{Y_i}\,e^{ - \sum_{j=1}^{i}Y_j}\right)^{p_2} \right]\,\PP(\tau_i+ \theta \in ds) \\[2mm] 
&&\leq \sum_{i=1}^{\infty} i^{\gamma\,q^*}\,\left\{\E\left[\left(e^{-\widehat{Y}}\right)^{p_2} \right]\right\}^{i-1}\,\int_0^t\,e^{s\,\phi_B(p_2)}\,\PP(\tau_i+ \theta \in ds) \\[2mm] 
&&\leq C^*\,\sum_{i=1}^{\infty} i^{\gamma\,q^*}\,\left\{\E\left[\left(e^{-\widehat{Y}}\right)^{p_2} \right]\right\}^{i-1}\,\PP(\tau_i+ \theta \leq t) \leq \dfrac{C^*}{\gamma\,q^*+1}\,\E\left[\left(N_{\Theta}(t)\right)^{\gamma\,q^* +1} \right] < \infty\,.
\eeao 

If $J_{F_A}^+ \geq 1$, for the application of \cite[Th. 2]{yi:chen:su:2011}, should hold
\beam \label{eq.KPY.4.25} 
\sum_{i=1}^{\infty} \left( \E\left[ \Theta_i^{p_1}(t)\,{\bf 1}_{\{\theta=0\}} \vee \Theta_i^{p_2}(t)\,{\bf 1}_{\{\theta=0\}}\right] \right)^{1/p_2} < \infty\,,
\eeam 
for all $t \in  \Lambda_T$. It is enough to show that the corresponding sum converges for each of the two terms. For the $p_2$ we obtain
\beam \label{eq.KPY.4.26} \notag
&&\sum_{i=1}^{\infty} \left( \E\left[ \Theta_i^{p_2}(t)\,{\bf 1}_{\{\theta=0\}} \right] \right)^{1/p_2} = \sum_{i=1}^{\infty} \left(\PP((\theta=0)\, \E\left[ \left( e^{-B(\tau_i)-\sum_{j=1}^{i-1} Y_j}\,{\bf 1}_{\{\tau_i \leq t\}} \right)^{p_2}\right] \right)^{1/p_2} \\[2mm] \notag
&&\leq \sum_{i=1}^{\infty} \left( \int_0^t \E\left[ (e^{-B(s)})^{p_2}\right]\,\E\left[ \left( e^{-\sum_{j=1}^{i-1} Y_j} \right)^{p_2}\right] \PP(\tau_i \in ds) \right)^{1/p_2}  \\[2mm] \notag 
&&\leq \sum_{i=1}^{\infty} \left( \left\{\E\left[ \left( e^{-\widehat{Y}} \right)^{p_2}\right] \right\}^{i-1}\,\int_0^t e^{-s\,\phi_B(p_2)}\,\PP(\tau_i \in ds) \right)^{1/p_2}\leq C^* \sum_{i=1}^{\infty} \left(\PP(\tau_i \leq t) \right)^{1/p_2}  \\[2mm]
&&= C^*\,\sum_{i=1}^{\infty} \left(\PP(N(t) \geq i) \right)^{1/p_2} \leq C^* \,\sum_{i=1}^{\infty} \left(\E\left[N^{\gamma\,q^*}(t) \right]\,
i^{-\gamma\,q^*} \right)^{1/p_2}  \\[2mm] \notag
&& \leq C^* \,\E\left[N^{\gamma\,q^*+1}(T) \right] \sum_{i=1}^{\infty}\,i^{-\gamma\,q^*/p_2} < \infty\,,
\eeam 
where $C^* \in (0,\,\infty)$ comes from $\phi_B(\cdot)$, at the sixth step we use Markov inequality. At the last step we recall that 
$p_2 > 1$, $q^* \geq p_2$. So $\gamma\,q^*/p_2 \geq \gamma$, and hence the last sum converges, since 
\beao
\sum_{i=1}^{\infty} \,i^{-\gamma\,q^*/p_2} \leq \sum_{i=1}^{\infty} \,i^{-\gamma}\,,
\eeao
because the last sum is a Cauchy series. By similar way we can work for $p_1$.

Following similar steps as with \eqref{eq.KPY.4.26}, we can show that
\beao
\sum_{i=1}^{\infty} \left( \E\left[ \widetilde{\Theta}_i^{q_1}(t)\,{\bf 1}_{\{\theta>0\}} \vee \widetilde{\Theta}_i^{q_2}(t)\,{\bf 1}_{\{\theta>0\}}\right] \right)^{1/p_2} < \infty\,,
\eeao
if $J_{G_A^{(i)}} \geq 1$.
 ~\halmos

The following lemma is useful for the proof of Corollary \ref{cor.KPY.3.1}.

\ble \label{lem.KPY.4.4}
Under the conditions of Corollary \ref{cor.KPY.3.1}, it holds 
\beam \label{eq.KPY.4.26.4} \notag
\lim \sup_{n \in \bbn} \left|\dfrac {\PP\left(\sum_{i=1}^n {\bf X}^{(i)}\,\Theta_i(t) \in x\,A \right)}{\mu(A) \bV(x) \sum_{i=1}^{n} \E\left[\Theta_i^{\alpha}(t)\,{\bf 1}_{\{\theta =0\}} \right] + \sum_{i=1}^{n} \nu_i(A) \overline{Q}_i(x) \E\left[\widetilde{\Theta}_i^{\beta}(t)\,{\bf 1}_{\{\theta >0\}} \right]} -1 \right| =0\,,\\
\eeam
uniformly for $t \in \Lambda_T$, for any fixed $T \in \Lambda$.
\ele

\pr~
At first, since $F \in MRV(\alpha,\,V,\,\mu)$, $G_i \in MRV(\beta,\,Q_i,\,\nu_i)$, with $i \in \bbn$ is implied that $F_A \in \mathcal{R}_{-\alpha}$, $G_A^{(i)} \in \mathcal{R}_{-\beta}$, with $i \in \bbn$. Thus we obtain 
\beao
J_{F_A}^- =J_{F_A}^+ = \alpha\,,
\eeao 
and 
\beao
J_{G_A^{(i)}}^- =J_{G_A^{(i)}}^+ = \beta\,.
\eeao 
From the moment condition \eqref{eq.KPY.3.8}, and through the monotone convergence theorem we find that for any $i \in \bbn$ it holds  
\beam \label{eq.KPY.4.26.5} \notag
&&\lim_{\delta \downarrow 0} \left(M_{\alpha - \delta}^{(i)} \vee M_{\alpha + \delta} ^{(i)}\right) = \lim_{\delta \downarrow 0} \int_0^{\infty} \left(s^{\alpha -\delta} \vee s^{\alpha + \delta}\right)\,\PP(e^{-Y_i} \in ds) \\[2mm]
&& =\int_0^{\infty} \lim_{\delta \downarrow 0} \left(s^{\alpha -\delta} \vee s^{\alpha + \delta}\right)\,\PP(e^{-Y_i} \in ds) =M_{\alpha} \leq 1\,.
\eeam
Quite similarly we obtain
\beam \label{eq.KPY.4.26.6} 
\lim_{\delta \downarrow 0} \left(M_{\beta - \delta}^{(i)} \vee M_{\beta + \delta} ^{(i)}\right) \leq 1\,,
\eeam
for any $i \in \bbn$. Hence, from \eqref{eq.KPY.4.26.5} and \eqref{eq.KPY.4.26.6}, there exist $\gamma,\,\gamma' >0$, such that 
\beao
\left(M_{\alpha - \gamma}^{(i)} \vee M_{\alpha + \gamma} ^{(i)}\right)\bigvee  \left(M_{\beta - \gamma'}^{(i)} \vee M_{\beta + \gamma'} ^{(i)}\right) \leq 1\,,
\eeao
for any $i \in \bbn$, that implies relation \eqref{eq.KPY.3.1} is satisfied, hence from Lemma \ref{lem.KPY.4.3} we obtain that \eqref{eq.KPY.4.20} is true.

Further, since the ${\bf X}^{(i)}$ and $\Theta_i^{\alpha}(t)\,{\bf 1}_{\{\theta =0\}}$ are independent, and through the first relation in \eqref{eq.KPY.4.5}, we obtain through Breiman's theorem (see, \cite{cline:samorodnitsky:1994}) that for any $i \in \bbn$, and for all $t \in \Lambda_T$, it holds
\beam \label{eq.KPY.4.26.7} \notag
&&\PP\left( {\bf X}^{(i)}\,\Theta_i(t)\,{\bf 1}_{\{\theta =0\}}  \in x\,A \right)=\PP\left( X_A^{(i)}\,\Theta_i(t)\,{\bf 1}_{\{\theta =0\}} >x \right) \\[2mm]
&&\sim \PP\left( X_A^{(i)} >x \right)\,\E[\Theta_i^{\alpha}(t)\,{\bf 1}_{\{\theta =0\}}] \sim  \mu(A) \bV(x) \E\left[\Theta_i^{\alpha}(t)\,{\bf 1}_{\{\theta =0\}} \right]\,.
\eeam 
Similarly, using the second relation in  \eqref{eq.KPY.4.5} we find that for any $i \in \bbn$, and for all $t \in \Lambda_T$, it holds 
\beam \label{eq.KPY.4.26.8} 
\PP\left( \widehat{\bf X}^{(i)}\,\Theta_i(t)\,{\bf 1}_{\{\theta =0\}}  \in x\,A \right) \sim  \nu_i(A) \overline{Q}_i(x) \E\left[\widetilde{\Theta}^{\beta}_i(t)\,{\bf 1}_{\{\theta =0\}} \right]\,.
\eeam

Let some fixed $N \in \bbn$. Then by Lemma \ref{lem.KPY.4.3}, and relations \eqref{eq.KPY.4.26.7} and  \eqref{eq.KPY.4.26.8} we obtain that it holds
\beam \label{eq.KPY.4.26.9} \notag
\lim \sup_{N \geq n \in \bbn} \left|\dfrac {\PP\left(\sum_{i=1}^n {\bf X}^{(i)}\,\Theta_i(t) \in x\,A \right)}{\mu(A) \bV(x) \sum_{i=1}^{n} \E\left[\Theta_i^{\alpha}(t) {\bf 1}_{\{\theta =0\}} \right] + \sum_{i=1}^{n} \nu_i(A) \overline{Q}_i(x) \E\left[\widetilde{\Theta}_i^{\beta}(t) {\bf 1}_{\{\theta >0\}} \right]} -1 \right| =0 ,\\
\eeam
 uniformly for $t \in \Lambda_T$.

It remains to show that \eqref{eq.KPY.4.26.4} holds also for all $n > N$. 
For the upper bound, in a similar way as for derivation of \eqref{eq.KPY.4.23},
 we obtain that for any $\vep >0$, there exists some $N \in \bbn$ large enough, such that for all $n > N$, $t \in \Lambda_T$ it holds
\beam \label{eq.KPY.4.26.10} \notag
&&\PP\left(\sum_{i=1}^n {\bf X}^{(i)}\,\Theta_i(t) \in x\,A \right)
 \\[2mm] \notag
&&\lesssim \sum_{i=1}^{\infty} \PP\left( {\bf X}^{(i)}\,\Theta_i(t)\,{\bf 1}_{\{\theta =0\}} \in x\,A \right) + \sum_{i=1}^{\infty} \PP\left( \widehat{\bf X}^{(i)}\,\widetilde{\Theta}_i(t)\,{\bf 1}_{\{\theta >0\}} \in x\,A \right)\\[2mm]
&& = \left( \sum_{i=1}^{N} +\sum_{i=N+1}^{\infty} \right)\PP\left( {\bf X}^{(i)}\,\Theta_i(t)\,{\bf 1}_{\{\theta =0\}} \in x\,A \right) + \left( \sum_{i=1}^{N} +\sum_{i=N+1}^{\infty} \right)  \\[2mm] \notag
&&\PP\left(\widehat{\bf X}^{(i)}\,\widetilde{\Theta}_i(t)\,{\bf 1}_{\{\theta >0\}} \in x\,A \right)=:I_{11}(x,\,t) +I_{12}(x,\,t) +I_{21}(x,\,t) +I_{22}(x,\,t)\,,
\eeam 
where at the first step we used similar techniques as in \eqref{eq.KPY.4.23}, through application of \cite[Th. 2]{yi:chen:su:2011}. 

For $I_{11}(x,\,t)$ and $I_{21}(x,\,t)$, applying \eqref{eq.KPY.4.26.7}, \eqref{eq.KPY.4.26.8}, for all $n > N$, $t \in \Lambda_T$ we obtain
\beam \label{eq.KPY.4.26.a} \notag
&&I_{11}(x,\,t) \sim \mu(A)\,\overline{V}(x) \sum_{i=1}^N \E\left[\Theta_i^{\alpha}(t)\,{\bf 1}_{\{\theta =0\}} \right] \leq \mu(A)\,\overline{V}(x) \sum_{i=1}^n \E\left[\Theta_i^{\alpha}(t)\,{\bf 1}_{\{\theta =0\}} \right]\,,\\[2mm]
&& 
I_{21}(x,\,t) \lesssim  \sum_{i=1}^n \nu_i(A)\,\overline{Q}_i(x)\, \E\left[\widetilde{\Theta}_i^{\beta}(t)\,{\bf 1}_{\{\theta >0\}} \right] \,.
\eeam

For $I_{12}(x,\,t)$ and $I_{22}(x,\,t)$, by similar way as in derivation of \eqref{eq.KPY.4.16}, \eqref{eq.KPY.4.19}, for any $\vep > 0$, there exists some $n_0 \in \bbn$, such that for all $N\geq n_0$, we find
\beam \label{eq.KPY.4.26.b} \notag
&&I_{12}(x,\,t) \lesssim \vep\,\PP({\bf X}^{(1)}\in x\,A ) \lesssim \vep\,C''\,\mu(A)\,\overline{V}(x) \,\E\left[\Theta_1^{\alpha}(t)\,{\bf 1}_{\{\theta =0\}} \right]
\,,\\[2mm]
&& 
I_{22}(x,\,t) \lesssim \vep\,C'' \,\nu_1(A)\,\overline{Q}_1(x)\, \E\left[\widetilde{\Theta}_1^{\beta}(t)\,{\bf 1}_{\{\theta >0\}} \right] \,.
\eeam
uniformly for $t \in \Lambda_T$, for some constant $C'' \in (0,\,\infty)$, that appears in combination with relations  \eqref{eq.KPY.4.26.7} and \eqref{eq.KPY.4.26.8}.

From \eqref{eq.KPY.4.26.10}, in combination with  \eqref{eq.KPY.4.26.a} and \eqref{eq.KPY.4.26.b}, for all $n > N$, $t \in \Lambda_T$ we obtain
\beam \label{eq.KPY.4.26.c}
&&\PP\left( \sum_{i=1}^N {\bf X}^{(i)}\,\Theta_i(t) \in x\,A \right) \\[2mm] \notag
&& \lesssim (1+\vep\,C'')\,\left(\mu(A)\,\overline{V}(x) \sum_{i=1}^n \E\left[\Theta_1^{\alpha}(t)\,{\bf 1}_{\{\theta =0\}} \right]+   \sum_{i=1}^n \nu_i(A)\,\overline{Q}_i(x)\, \E\left[\widetilde{\Theta}_i^{\beta}(t)\,{\bf 1}_{\{\theta >0\}} \right] \right)\,.
\eeam
We note that the convergence of the last sum is implied by the fact that the $\nu_i$ are Radon measures, in combination with relation \eqref{eq.KPY.4.18}.

For the lower bound, via a similar way as for derivation of \eqref{eq.KPY.4.21}, with the help of \eqref{eq.KPY.4.26.7} and \eqref{eq.KPY.4.26.8}, we obtain that for all $n >  N$ it holds
\beam \label{eq.KPY.4.26.11} \notag
&&\PP\left(\sum_{i=1}^n {\bf X}^{(i)}\,\Theta_i(t) \in x\,A \right)
\gtrsim \mu(A) \bV(x) \sum_{i=1}^{N} \E\left[\Theta_i^{\alpha}(t)\,{\bf 1}_{\{\theta =0\}} \right]  \\[2mm] \notag
&&+ \sum_{i=1}^{N} \nu_i(A) \overline{Q}_i(x) \E\left[\widetilde{\Theta}_i^{\beta}(t)\,{\bf 1}_{\{\theta >0\}} \right] \geq \mu(A) \bV(x) \left(\sum_{i=1}^{n}- \sum_{i=N+1}^{\infty}\right)\E\left[\Theta_i^{\alpha}(t)\,{\bf 1}_{\{\theta =0\}} \right] \\[2mm]
&& + \left( \sum_{i=1}^{n} -\sum_{i=N+1}^{\infty} \right) \nu_i(A) \overline{Q}_i(x) \E\left[\widetilde{\Theta}_i^{\beta}(t)\,{\bf 1}_{\{\theta >0\}} \right]\\[2mm] \notag
&& \geq  (1-\vep) \left( \mu(A) \bV(x) \sum_{i=1}^{n} \E\left[\Theta_i^{\alpha}(t)\,{\bf 1}_{\{\theta =0\}} \right] + \sum_{i=1}^{n} \nu_i(A) \overline{Q}_i(x) \E\left[\widetilde{\Theta}_i^{\beta}(t)\,{\bf 1}_{\{\theta >0\}} \right] \right)\,,
\eeam 
uniformly for $t \in  \Lambda_T$, where at the last step we used the \eqref{eq.KPY.4.14} and \eqref{eq.KPY.4.18} and the fact that $\nu_i$ are Radon measures. From \eqref{eq.KPY.4.26.c} and  \eqref{eq.KPY.4.26.11} and the arbitrary choice of $\vep>0$ we find 
\beam \label{eq.KPY.4.26.12} \notag
\lim \sup_{N < n \in \bbn} \left|\dfrac {\PP\left(\sum_{i=1}^n {\bf X}^{(i)}\,\Theta_i(t) \in x\,A \right)}{\mu(A) \bV(x) \sum_{i=1}^{n} \E\left[\Theta_i^{\alpha}(t) {\bf 1}_{\{\theta =0\}} \right] + \sum_{i=1}^{n} \nu_i(A) \overline{Q}_i(x) \E\left[\widetilde{\Theta}_i^{\beta}(t) {\bf 1}_{\{\theta >0\}} \right]} -1 \right| =0 ,\\
\eeam
uniformly for $t \in  \Lambda_T$. From \eqref{eq.KPY.4.26.9} and  \eqref{eq.KPY.4.26.12} we conclude that \eqref{eq.KPY.4.26.4} is true, uniformly for $t \in  \Lambda_T$. 
~\halmos

\noindent{\bf Proof of Theorem \ref{th.KPY.3.1}.}~
Applying Lemma \ref{lem.KPY.4.3} for $n=\infty$, we obtain
\beam \label{eq.KPY.4.27} \notag
&&\PP({\bf D}(t) \in x\,A) = \PP \left(\sum_{i=1}^{\infty} {\bf X}^{(i)}\,e^{-B(\tau_i+\theta)-J[(\tau_i + \theta)-]}\,{\bf 1}_{\{\tau_i+\theta \leq t\}} \in x\,A \right) \\[2mm] \notag
&&=\PP\left( \sum_{i=1}^{\infty}  {\bf X}^{(i)}\,\Theta_i(t) \in x\,A 
 \right) \sim \sum_{i=1}^{\infty}  \PP\left( {\bf X}^{(i)}\,\Theta_i(t)\,{\bf 1}_{\{\theta=0\}} \in x\,A  \right)\\[2mm] \notag
&&+\sum_{i=1}^{\infty}  \PP\left( \widehat{\bf X}^{(i)}\,\widetilde{\Theta}_i(t)\,{\bf 1}_{\{\theta>0\}} \in x\,A  \right) = \sum_{i=1}^{\infty} \PP({\bf X}^{(i)}\,\Theta_i(t) \in x\,A) \,, 
\eeam 
uniformly for $t \in  \Lambda_T$. From \eqref{eq.KPY.4.27} we get the desired result.
~\halmos

\noindent{\bf Proof of Corollary \ref{cor.KPY.3.1}.}~
It is implied, similarly with the proof of Theorem \ref{th.KPY.3.1}, but by application of Lemma \ref{lem.KPY.4.4}, instead of Lemma \ref{lem.KPY.4.3}. 
~\halmos

\section{Dependence between insurance and financial risks} \label{sec.KPT29U.5}

One of the interesting questions on the conditions of Theorem  \ref{th.KPY.3.1}, is when holds the condition  $G_i \in (\mathcal{C} \cap \mathcal{P_D})_A$. In this section we present sufficient conditions for these inclusions, employing closure properties with respect to scalar product. Further, under a weak, but general, dependence structure, we show that relation \eqref{eq.KPY.3.7} can be simplified by reduction of the $\widehat{\bf X}^{(i)} ={\bf X}^{(i)}\,e^{-Y_i}$ to product of independent random vectors using solution of the dependence.

Finally, we discuss two cases, where the conditions of Corollary  \ref{cor.KPY.3.1} are satisfied, under a weak and a strong dependence structures, depicting  directly the essential influence of dependence between the two risks on the asymptotic estimations. Hereafter, we denote for sake of convenience, $Z_i:=e^{-Y_i}$ and $Z_i \stackrel{d}{\sim} H_i$, for any $i \in \bbn$.

\subsection{Weak dependence on class $(\mathcal{C} \cap \mathcal{P_D})_A$} \label{subsec.KPT29U.5.1}

The following dependence structure appeared in a one-dimensional framework through copulas in \cite{asimit:jones:2008}, and further was used and extended to several applications on discrete and continuous time risk models, see for example \cite{li:tang:wu:2010}, \cite{tang:yuan:2016}, \cite{yuan:lu:2023}, \cite{chen:cheng:2024} among many others. We shall use here a corresponding multidimensional projection of this concrete dependence structure for the modeling of the dependence between the insurance claims ${\bf X}^{(i)}$ and the financial jumps of the investment portfolio $Z_i=e^{-Y_i}$.

\begin{assumption} \label{ass.KPY.5.1}
Let consider some fixed set  $A \in \mathscr{R}$. For any $i \in \bbn$, we suppose that there exists some measurable function $h_i\;:\;[0,\,\infty) \to (0,\,\infty)
$, such that it holds
\beam \label{eq.KPY.5.1}  
\inf_{z \in E} h_i(z) >0\,,
\eeam 
where $E$ a neighborhood of the right endpoint of $S(Z_i)$, such that we obtain
\beam \label{eq.KPY.5.2}  
\PP({\bf X}^{(i)} \in x\,A\;|\; Z_i=z) \sim h_i(z)\,\PP({\bf X}^{(i)} \in x\,A)\,,
\eeam 
uniformly for $z \in S(Z_i)$.   
\end{assumption}

\bre \label{rem.KPY.5.1}
From \cite[Prop. 2.4]{cui:wang:2024} we find that the function $h_i$, for $i \in \bbn$, is bounded from above, that means that there exists some constant $C>0$ such that for any $i \in \bbn$, $z \in S(Z_i)$ we obtain 
\beam \label{eq.KPY.5.3} 
h_i(z) \leq C\,.
\eeam
Furthermore, by integration on both sides of \eqref{eq.KPY.5.1}, with respect to $\PP(Z_i \in dz)$, over the $S(Z_i)$ we find $\E[h_i(Z_i)]=1$.
\ere

The following result provides the sufficient conditions for the inclusion $G_i \in (\mathcal{C} \cap \mathcal{P_D})_A$. We apply Assumption \ref{ass.KPY.5.1} supposing that for the pair $({\bf X},\,Z)$ are satisfied relations \eqref{eq.KPY.5.1} and \eqref{eq.KPY.5.2}, for some measurable function $h(\cdot)\;:\;[0,\,\infty) \to (0,\,\infty)$. This result can be derived as special case of \cite[Th. 4.1]{konstantinides:passalidis:2025h}, taking into account the \cite[Cor. 2.1]{tang:2006a}. However, here we follow an alternative proof, that will be partially useful in the proof of Corollary \ref{cor.KPY.5.1}.

\bpr \label{pr.KPY.5.1}
Let $({\bf X},\,Z)$ satisfies the Assumption \ref{ass.KPY.5.1}. If $F \in (\mathcal{C} \cap \mathcal{P_D})_A$ and $\E[Z^p] < \infty$ for some $p > J_{F_A}^+$, then we obtain
$\widehat{\bf X} := {\bf X}\,Z  \stackrel{d}{\sim} G \in (\mathcal{C} \cap \mathcal{P_D})_A$.  
\epr

For the proof of the previous proposition (see, in Subsection \ref{subsec.KPT29U.5.2}), as also for the formulation of the following corollary we need the new random variables $\{Z_i^{*}=e^{-Y_i^*}\,,\; i \in \bbn\}$, that are independent of $\{\widehat{\bf X}^{(i)}\,,\; i \in \bbn \}$ and $Z_i^{*} \stackrel{d}{\sim} H_i^{*}$, where 
\beam \label{eq.KPY.5.4}  
H_i^{*}(dz) = h_i(z)\,H_i(dz)\,.
\eeam 

From Remark \ref{rem.KPY.5.1} we can see that the $\{H_i^{*}\,,\;i \in \bbn\}$ are proper distributions, so from \eqref{eq.KPY.5.1} and \eqref{eq.KPY.5.3} follows that the $(H_i^{*}\,,\;H_i)$ for $i \in \bbn$, are equivalent measures, namely for any $i \in \bbn$ and some Borel mearure $\bbb$ it holds
\beao
\PP(Z_i^* \in \bbb)= 0 \;\Longleftrightarrow\; \PP(Z_i \in \bbb)= 0\,.
\eeao 

In Theorem \ref{th.KPY.3.1}, see also Remark \ref{rem.KPY.3.4}, relation \eqref{eq.KPY.3.7} contains independent products of random variable, however it also contains and the $\widehat{\bf X}^{(i)}= {\bf X}^{(i)}\,e^{-Y_i}$. Unfortunately, in some cases the calculation of the distribution of $\widehat{\bf X}^{(i)}$, due to the fact that it is product of dependent random vectors, becomes difficult. The next result, provides an elegant form, in which are included the products of independent random vectors in each probability.  

\bco  \label{cor.KPY.5.1}
 Under the conditions of Theorem \ref{th.KPY.3.1}, with the only restriction that  the $\{({\bf X}^{(i)},\,e^{-Y_i})\,,\;i \in \bbn\}$ satisfy Assumption  \ref{ass.KPY.5.1}, it holds
\beam \label{eq.KPY.5.5}  
\PP({\bf D}(t) \in x\,A) \sim\sum_{i=1}^{\infty} \PP({\bf X}^{(i)} \Theta_i(t) {\bf 1}_{\{\theta =0\}} \in x\,A) +  \sum_{i=1}^{\infty} \PP({\bf X}^{(i)} e^{-Y_i^*}\widetilde{\Theta}_i(t) {\bf 1}_{\{\theta>0\}} \in x\,A) ,
\eeam 
uniformly for $t \in \Lambda_T$, for any fixed $T \in \Lambda$.
\eco

\subsection{Argumentation of subsection \ref{subsec.KPT29U.5.1}} \label{subsec.KPT29U.5.2}

In this subsection we provide the proof of Proposition \ref{pr.KPY.5.1} and of Corollary \ref{cor.KPY.5.1}. We start with a preliminary lemma, that has its own right, since it provides the solution of the dependence from Assumption \ref{ass.KPY.5.1} in the class $\mathcal{D}_A$, for $F$, under moment conditions for the $Z$. We recall that $Z^* \stackrel{d}{\sim} H^*$, with $H^*(dz) = h(z)\,H(dz)$, 
while the $Z^*$ is independent of ${\bf X}$.

\ble \label{lem.KPY.5.1}
Let $({\bf X},\,Z)$ satisfies Assumption \ref{ass.KPY.5.1}, with $F \in \mathcal{D}_A$, and for some $p > J_{F_A}^+$ it holds $\E[Z^p] < \infty$. Then it holds
\beam \label{eq.KPY.5.6}  
\PP({\bf X}\,Z \in x\,A) \sim \PP({\bf X}\,Z^* \in x\,A) \,.
\eeam 
\ele 

\pr~
Without loss of generality, we assume that $S(Z)=[0,\,\infty)$, as in case of $S(Z) \subsetneq [0,\,\infty)$ the same argumentation remains true (and it is easier). Let choose some $k \in (J_{F_A}^+ /p ,\,1)$. Then  it holds
\beam \label{eq.KPY.5.7}  
&&\PP({\bf X}\,Z \in x\,A) = \int_0^{\infty} \PP\left({\bf X} \in \dfrac xz \,A\;\big|\;Z =z\right)\,\PP(Z \in dz) \\[2mm] \notag
&&=\left( \int_0^{x^k} + \int_{x^k}^{\infty} \right)  \PP\left(X_A > \dfrac xz\;\big|\;Z =z\right)\,\PP(Z \in dz)=:I_1(x,\,k)+ I_2(x,\,k)\,.
\eeam 
For $I_2(x,\,k)$, via application of Markov inequality we obtain
\beam \label{eq.KPY.5.8}  
I_2(x,\,k) \leq \PP(Z > x^k) \leq \E[Z^p]\,x^{-k\,p} =o\left[\bF_A(x) \right] =o\left[\PP({\bf X} \in x\,A) \right]\,,
\eeam 
where at the third step we took into consideration \eqref{eq.KPY.2.10}. Next, since $\E[Z^p] < \infty$, from \eqref{eq.KPY.5.4} and \eqref{eq.KPY.5.3} we find
\beam \label{eq.KPY.5.9}  
\E\left[(Z^*)^p \right] =\int_0^{\infty}z^p\,H^*(dz) =\int_0^{\infty}z^p\,h(z)\,H(dz) \leq C\,\E\left[Z^p\right] < \infty\,.
\eeam
So, taking into account that $F \in \mathcal{D}_A$ and $Z^*$ is independent of ${\bf X}$, it follows by \cite[Th. 3.3(iv)]{cline:samorodnitsky:1994} that
\beam \label{eq.KPY.5.10}  
\PP\left({\bf X}\,Z^* \in x\,A \right) =\PP\left(X_A\,Z^* >x \right) \asymp \PP\left(X_A >x \right)= \PP\left({\bf X} \in x\,A \right) \,.
\eeam
From \eqref{eq.KPY.5.8} and \eqref{eq.KPY.5.10} we obtain
\beam \label{eq.KPY.5.11}  
I_2(x,\,k)=\left[\PP\left({\bf X}\,Z^* \in x\,A  \right) \right] \,.
\eeam

For $I_1(x,\,k)$, we have
\beam \label{eq.KPY.5.12}  
&&I_1(x,\,k) \sim \int_0^{x^k} h(z)\PP\left( X_A >\dfrac xz \right)\,H(dz) \\[2mm] \notag
&&=\left( \int_0^{\infty} - \int_{x^k}^{\infty} \right)  \PP\left( X_A > \dfrac xz \right)\,H^*(dz) \\[2mm] \notag
&&=(1-o(1)) \int_0^{\infty} \PP\left( X_A > \dfrac xz \right)\,H^*(dz) =(1-o(1))  \PP\left( {\bf X}\,Z^* \in x\,A  \right)\,,
\eeam 
where at the first step  we applied the dominated convergence theorem, while at the third step we used the relation
\beao
J_2^*=\int_{x^k}^{\infty} \PP\left( X_A > \dfrac xz \right)\,H^*(dz)=o\left[  \PP\left( {\bf X} \in x\,A  \right) \right]=o\left[  \PP\left( {\bf X}\,Z^* \in x\,A  \right) \right]\,,
\eeao
that is implies following similar path with the derivation of relation \eqref{eq.KPY.5.11}. From \eqref{eq.KPY.5.11} and \eqref{eq.KPY.5.12}, in combination with \eqref{eq.KPY.5.7}, we reach to relation  \eqref{eq.KPY.5.6}. 
~\halmos

Now we can proceed to the proof of Proposition  \ref{pr.KPY.5.1}.

\noindent{\bf Proof of Proposition  \ref{pr.KPY.5.1}.}~
From Lemma \ref{lem.KPY.5.1} we obtain
\beam \label{eq.KPY.5.13}  
\PP\left( X_A\,Z >  x \right)\sim \PP\left( X_A\,Z^* > x \right)\,.
\eeam 

Since $F_A \in \mathcal{C} \cap \mathcal{P_D}$, by \cite[Th. 3.4(ii)]{cline:samorodnitsky:1994} and \cite[Th. 5.1(i)]{konstantinides:passalidis:2024d} we find that the distribution of $X_A\,Z^*$ belongs to class 
$ \mathcal{C}\cap \mathcal{P_D}$. Hence, from \eqref{eq.KPY.5.13} and the closure property of $ \mathcal{C}\cap \mathcal{P_D}$ with respect to strong tail equivalence, we obtain that 
$X_A\,Z \stackrel{d}{\sim} G_A \in \mathcal{C}\cap \mathcal{P_D}$, so we get 
$G \in \mathcal{C}\cap \mathcal{P_D}$.
~\halmos

Before we proceed to the proof of Corollary \ref{cor.KPY.5.1} we need an auxiliary lemma.

\ble \label{lem.KPY.5.2}
Under the conditions of Corollary \ref{cor.KPY.5.1}, for any $i \in\bbn$ it holds 
\beam \label{eq.KPY.5.14}  
\PP\left(\widehat{\bf X}^{(i)}\,\widetilde{\Theta}_i(t)\,{\bf 1}_{\{\theta>0\}} \in  x\,A \right) \sim \PP\left({\bf X}^{(i)}\,e^{-Y_i^*}\,\widetilde{\Theta}_i(t)\,{\bf 1}_{\{\theta>0\}} \in  x\,A \right)\,,
\eeam
uniformly for $t \in \Lambda_T$, for any fixed $T \in \Lambda$.
\ele

\pr~
At first, for all $i \in \bbn$, $t \in \Lambda_T$, and some $k \in (J_{F_A}^+ /p_2 ,\,1)$ it holds
\beam \label{eq.KPY.5.15}  
&&\PP\left(\widehat{\bf X}^{(i)}\,\widetilde{\Theta}_i(t)\,{\bf 1}_{\{\theta>0\}} \in  x\,A \right) = \PP\left( X_A^{(i)}\,e^{-Y_i}\,\widetilde{\Theta}_i(t)\,{\bf 1}_{\{\theta>0\}} > x \right)  \\[2mm] \notag
&&=\left(\int_0^{x^k} + \int_{x^k}^{\infty} \right) \PP\left( X_A^{(i)}\,\widetilde{\Theta}_i(t)\,{\bf 1}_{\{\theta>0\}} >\dfrac xz \;\bigg|\;e^{-Y_i}=z\right)\,\PP\left( e^{-Y_i} \in dz \right) \\[2mm] \notag
&&= J_1(x,\,t;\,k) +J_2(x,\,t;\,k)\,.
\eeam 
Let us note that since for all $t \in \Lambda_T$, $\E\left[\widetilde{\Theta}_i^{p_2}(t)\,{\bf 1}_{\{\theta>0\}} \right] < \infty$, (see, relation \eqref{eq.KPY.4.5}), and $\widetilde{\Theta}_i(t)$ is independent of $X_A \stackrel{d}{\sim} F_A \in \mathcal{D} $, from 
\cite[Lem. 3.9]{tang:tsitsiashvili:2003} follows that the product $X_A^{(i)}\,\widetilde{\Theta}_i(t)\,{\bf 1}_{\{\theta>0\}}$ has upper Matuszewska index equal to $J_{F_A}^+$. Further, for the $J_2(x,\,t;\,k)$ and for all $t \in \Lambda_T$ it holds
\beam \label{eq.KPY.5.16}  
&&J_2(x,\,t;\,k) \leq \PP\left(e^{-Y_i} > x^k \right) \leq \E\left[ (e^{-Y_i})^{p_2} \right]\,x^{-k\,p_2} = M_{p_2}^{(i)}\,x^{-k\,p_2} \\[2mm] \notag
&&=o\left[ \PP\left( X_A^{(i)}\,\widetilde{\Theta}_i(t)\,{\bf 1}_{\{\theta>0\}} > x\right) \right] =o\left[ \PP\left( {\bf X}^{(i)}\,\widetilde{\Theta}_i(t)\,{\bf 1}_{\{\theta>0\}} \in x\,A \right) \right]\,,
\eeam
where at the second step we used Markov inequality and in the pre-last step we used \eqref{eq.KPY.2.10}. From the other hand side, since for any $i \in \bbn$, 
\beao
&& \E\left[ (e^{-Y_i^*})^{p_2} \right]= \int_0^{\infty} z^{p_2} \,H_i^*(dz) \\[2mm] \notag
&&=\int_0^{\infty} z^{p_2} \,h_i(z)\,H_i(dz)  \leq C\,M_{p_2}^{(i)} \leq C < \infty\,,
\eeao
where at the third step we used relation \eqref{eq.KPY.5.3}, and at the fourth step relation \eqref{eq.KPY.3.1}. From the previous relation and \cite[Th. 3.3(iv)]{cline:samorodnitsky:1994} for all $t \in \Lambda_T$ we obtain
\beam \label{eq.KPY.5.17}  
\PP\left( {\bf X}^{(i)}\,\widetilde{\Theta}_i(t)\,{\bf 1}_{\{\theta>0\}} \in x\,A \right) \asymp \PP\left( {\bf X}^{(i)}\,e^{-Y_i^*}\,\widetilde{\Theta}_i(t)\,{\bf 1}_{\{\theta>0\}} \in x\,A \right)\,.
\eeam
Therefore, from \eqref{eq.KPY.5.16} and \eqref{eq.KPY.5.17} for all $t \in \Lambda_T$ we find
\beam \label{eq.KPY.5.18}  
J_2(x,\,t;\,k)=o\left[ \PP\left( {\bf X}^{(i)}\,e^{-Y_i^*}\,\widetilde{\Theta}_i(t)\,{\bf 1}_{\{\theta>0\}} \in x\,A \right) \right]\,.
\eeam

Now we deal with $J_1(x,\,t;\,k)$. For all $t \in \Lambda_T$ we obtain
\beam \label{eq.KPY.5.19}   \notag
&&J_1(x,\,t;\,k) =\int_0^{x^k} \int_0^{\infty} \PP\left( X_A > \dfrac x{u\,z} \;\big|\; e^{-Y_i} =z\right) \,\PP\left(\widetilde{\Theta}_i(t)\,{\bf 1}_{\{\theta>0\}} \in du \right) \,\PP\left( e^{-Y_i} \in dz \right)\\[2mm] \notag
&&=\int_0^{x^k} \left(\int_0^{x^k} + \int_{x^k}^{\infty} \right) \PP\left( X_A > \dfrac x{u\,z} \;\big|\; e^{-Y_i} =z\right) \,\PP\left(\widetilde{\Theta}_i(t)\,{\bf 1}_{\{\theta>0\}} \in du \right) \,\PP\left( e^{-Y_i} \in dz \right)\\[2mm]
&&=J_{11}(x,\,t;\,k) + J_{12}(x,\,t;\,k)\,.
\eeam
For $J_{12}(x,\,t;\,k)$, for all $t \in \Lambda_T$, by Markov inequality we find
\beam \label{eq.KPY.5.20}   \notag
&&J_{12}(x,\,t;\,k) \leq \int_0^{x^k} \PP\left(\widetilde{\Theta}_i(t)\,{\bf 1}_{\{\theta>0\}} > x^k \right) \,\PP\left( e^{-Y_i} \in dz \right)\leq  \PP\left(\widetilde{\Theta}_i(t)\,{\bf 1}_{\{\theta>0\}} > x^k \right) \\[2mm]
&&\leq \E\left[ \left(\widetilde{\Theta}_i(t)\right)^{p_2}\,{\bf 1}_{\{\theta>0\}}\right]\,x^{-k\,p_2} =o\left[ \PP\left( {\bf X}^{(i)}\,\widetilde{\Theta}_i(t)\,{\bf 1}_{\{\theta>0\}} \in x\,A \right) \right] \\[2mm] \notag
&& =o\left[ \PP\left( {\bf X}^{(i)}\,e^{-Y_i^*}\,\widetilde{\Theta}_i(t)\,{\bf 1}_{\{\theta>0\}} \in x\,A \right) \right] \,,
\eeam
where at the last step we employed relation \eqref{eq.KPY.5.17}. 

From Assumption \ref{ass.KPY.5.1} and through dominated convergence theorem, for all $t \in \Lambda_T$ we obtain
\beam \label{eq.KPY.5.21}   \notag
&&J_{11}(x,\,t;\,k) \sim \int_0^{x^k} \int_0^{x^k} h_i(z)\,\PP\left( X_A^{(i)} > \dfrac x{u\,z} \right) \,\PP\left(\widetilde{\Theta}_i(t)\,{\bf 1}_{\{\theta>0\}} \in du \right)\,\PP\left( e^{-Y_i} \in dz \right)\\[2mm] \notag
&& =\int_0^{x^k} \left( \int_0^{\infty} - \int_{x^k}^{\infty} \right) h_i(z)\,\PP\left( X_A^{(i)} > \dfrac x{u\,z} \right) \,\PP\left(\widetilde{\Theta}_i(t)\,{\bf 1}_{\{\theta>0\}} \in du \right)\,\PP\left( e^{-Y_i} \in dz \right)\\[2mm]
&& = [1-o(1)]\left( \int_0^{\infty} - \int_{x^k}^{\infty} \right) h_i(z)\,\PP\left( X_A^{(i)}\,\widetilde{\Theta}_i(t)\,{\bf 1}_{\{\theta>0\}}  > \dfrac x{z} \right) \,\PP\left( e^{-Y_i} \in dz \right) \\[2mm] \notag
&& = [1-2\,o(1)] \int_0^{\infty}h_i(z)\,\PP\left({\bf X}^{(i)}\,\widetilde{\Theta}_i(t)\,{\bf 1}_{\{\theta>0\}}  > \dfrac x{z}\,A \right) \,\PP\left( e^{-Y_i} \in dz \right) \\[2mm] \notag
&& = [1-2\,o(1)]\,\PP\left({\bf X}^{(i)}\, e^{-Y_i^*}\,\widetilde{\Theta}_i(t)\,{\bf 1}_{\{\theta>0\}}  >  x\,A \right) \,,
\eeam
where at the third step, we used similar arguments with that for \eqref{eq.KPY.5.20}, in combination with \eqref{eq.KPY.5.3}. At the fourth step we follow the same path of  \eqref{eq.KPY.5.18}, in combination with  \eqref{eq.KPY.5.3}.

From  \eqref{eq.KPY.5.19} -  \eqref{eq.KPY.5.21}, for all $t \in \Lambda_T$ we find
\beam \label{eq.KPY.5.22}  
J_1(x,\,t;\,k) \sim \PP\left( {\bf X}^{(i)}\,e^{-Y_i^*}\,\widetilde{\Theta}_i(t)\,{\bf 1}_{\{\theta>0\}} \in x\,A \right) \,.
\eeam
Putting  \eqref{eq.KPY.5.18} and \eqref{eq.KPY.5.22}, in  \eqref{eq.KPY.5.15}
we have \eqref{eq.KPY.5.14} uniformly for $t \in \Lambda_T$.
~\halmos

\noindent{\bf Proof of Corollary  \ref{cor.KPY.5.1}.}~
From Proposition \ref{pr.KPY.5.1}, all the assumptions of Theorem \ref{th.KPY.3.1} are satisfied, hence relation \eqref{eq.KPY.3.7} holds. It remains to show that for all $t \in \Lambda_T$, it holds
\beam \label{eq.KPY.5.23}  
\sum_{i=1}^{\infty} \PP\left( \widehat{\bf X}^{(i)}\,\widetilde{\Theta}_i(t)\,{\bf 1}_{\{\theta>0\}} \in x\,A \right) \sim \sum_{i=1}^{\infty} \PP\left( {\bf X}^{(i)}\,e^{-Y_i^*}\,\widetilde{\Theta}_i(t)\,{\bf 1}_{\{\theta>0\}} \in x\,A \right) \,.
\eeam  
From  \eqref{eq.KPY.5.14}, follows that for any fixed $N \in \bbn$, it holds
\beam \label{eq.KPY.5.24}  
\sum_{i=1}^{n} \PP\left( \widehat{\bf X}^{(i)}\,\widetilde{\Theta}_i(t)\,{\bf 1}_{\{\theta>0\}} \in x\,A \right) \sim \sum_{i=1}^{n} \PP\left( {\bf X}^{(i)}\,e^{-Y_i^*}\,\widetilde{\Theta}_i(t)\,{\bf 1}_{\{\theta>0\}} \in x\,A \right) \,,
\eeam   
uniformly for $1\leq n \leq N$, and $t \in \Lambda_T$. It remains to show that \eqref{eq.KPY.5.24} holds uniformly for $n> N$, and $t \in \Lambda_T$.

For all $n> N$, and $t \in \Lambda_T$, it holds
\beam \label{eq.KPY.5.25}  
&&\sum_{i=1}^{n} \PP\left( \widehat{\bf X}^{(i)}\,\widetilde{\Theta}_i(t)\,{\bf 1}_{\{\theta>0\}} \in x\,A \right) \geq \sum_{i=1}^{N} \PP\left( \widehat{\bf X}^{(i)}\,\widetilde{\Theta}_i(t)\,{\bf 1}_{\{\theta>0\}} \in x\,A \right) \\[2mm] \notag
&& \sim \sum_{i=1}^{N} \PP\left( {\bf X}^{(i)}\,e^{-Y_i^*}\,\widetilde{\Theta}_i(t)\,{\bf 1}_{\{\theta>0\}} \in x\,A \right)  \\[2mm] \notag
&& = \left(\sum_{i=1}^{n} - \sum_{i=N+1}^n \right) \PP\left( {\bf X}^{(i)}\,e^{-Y_i^*}\,\widetilde{\Theta}_i(t)\,{\bf 1}_{\{\theta>0\}} \in x\,A \right) =: K_1(x,\,t;\,N) - K_2(x,\,t;\,N)\,,
\eeam 
where at the second step we used relation \eqref{eq.KPY.5.24} for $n=N$. From \eqref{eq.KPY.4.16}, for any $\vep > 0$, there exists some $N \in \bbn$, such that  for all $t \in \Lambda_T$, we obtain
\beam \label{eq.KPY.5.26}  
&&K_2(x,\,t;\,N) \leq \sum_{i=N+1}^{\infty} \PP\left( {\bf X}^{(i)}\,e^{-Y_i^*}\,\widetilde{\Theta}_i(t)\,{\bf 1}_{\{\theta>0\}} \in x\,A \right) \\[2mm] \notag
&& \lesssim \vep\,C^*\, \PP\left( \widehat{\bf X}^{(1)}\,\widetilde{\Theta}_1(t)\,{\bf 1}_{\{\theta>0\}} \in x\,A \right) \sim \vep\,C^*\, \PP\left( {\bf X}^{(1)}\,e^{-Y_1^*}\,\widetilde{\Theta}_1(t)\,{\bf 1}_{\{\theta>0\}} \in x\,A \right) \,,
\eeam
for some constant $C^*>0$. We note that the usage of \eqref{eq.KPY.4.19} is implied by Proposition \ref{pr.KPY.5.1}. From  \eqref{eq.KPY.5.25}, \eqref{eq.KPY.5.26} and the arbitrary choice of $\vep >0$, we see that for all $t \in \Lambda_T$, $n > N$ it holds
\beam \label{eq.KPY.5.27}  
&&\sum_{i=1}^{n} \PP\left( \widehat{\bf X}^{(i)}\,\widetilde{\Theta}_i(t)\,{\bf 1}_{\{\theta>0\}} \in x\,A \right) \gtrsim \sum_{i=1}^{n}  \PP\left( {\bf X}^{(i)}\,e^{-Y_i^*}\,\widetilde{\Theta}_i(t)\,{\bf 1}_{\{\theta>0\}} \in x\,A \right) \,,
\eeam  
From the other hand side, for all $t \in \Lambda_T$, $n > N$ we obtain
\beam \label{eq.KPY.5.28}   \notag
&&\sum_{i=1}^{n} \PP\left( \widehat{\bf X}^{(i)}\,\widetilde{\Theta}_i(t)\,{\bf 1}_{\{\theta>0\}} \in x\,A \right) \leq \left(\sum_{i=1}^{N} + \sum_{i=N+1}^{\infty} \right)  \PP\left(\widehat{\bf X}^{(i)}\,\widetilde{\Theta}_i(t)\,{\bf 1}_{\{\theta>0\}} \in x\,A \right)  \\[2mm]
&&=: L_1(x,\,t;\,N)+ L_2(x,\,t;\,N) \,.
\eeam
From relation  \eqref{eq.KPY.5.14}, for all $t \in \Lambda_T$, $n > N$ we find
\beam \label{eq.KPY.5.29}  
&&L_1(x,\,t;\,N) \sim \sum_{i=1}^{N} \PP\left( \widehat{\bf X}^{(i)}\,e^{-Y_i^*}\,\widetilde{\Theta}_i(t)\,{\bf 1}_{\{\theta>0\}} \in x\,A \right)  \\[2mm] \notag
&&\leq \sum_{i=1}^{n} \PP\left(\widehat{\bf X}^{(i)}\,e^{-Y_i^*}\,\widetilde{\Theta}_i(t)\,{\bf 1}_{\{\theta>0\}} \in x\,A \right) \,,
\eeam
and through similar way as with derivation of \eqref{eq.KPY.4.19}, we obtain that for any $\vep >0$, there exists $N \in \bbn$, such that for all $n >N$, and $t \in \Lambda_T$, it holds
\beam \label{eq.KPY.5.30}  
&&L_2(x,\,t;\,N) \leq \vep\,C^* \PP\left( \widehat{\bf X}^{(1)}\,\widetilde{\Theta}_1(t)\,{\bf 1}_{\{\theta>0\}} \in x\,A \right)  \\[2mm] \notag
&&\sim \vep\,C^* \PP\left({\bf X}^{(1)}\,e^{-Y_1^*}\,\widetilde{\Theta}_1(t)\,{\bf 1}_{\{\theta>0\}} \in x\,A \right) \,.
\eeam

From relations \eqref{eq.KPY.5.28} - \eqref{eq.KPY.5.30}, and the arbitrary choice of $\vep > 0$, we find that for all $n >N$, and $t \in \Lambda_T$, it holds
\beam \label{eq.KPY.5.31}  
\sum_{i=1}^{n} \PP\left( \widehat{\bf X}^{(i)}\,\widetilde{\Theta}_i(t)\,{\bf 1}_{\{\theta>0\}} \in x\,A \right) \lesssim \sum_{i=1}^{n} \PP\left({\bf X}^{(i)}\,e^{-Y_i^*}\,\widetilde{\Theta}_i(t)\,{\bf 1}_{\{\theta>0\}} \in x\,A \right) \,.
\eeam   

From relations \eqref{eq.KPY.5.24}, \eqref{eq.KPY.5.27}, \eqref{eq.KPY.5.31}, we have
\beam \label{eq.KPY.5.32}  
\sum_{i=1}^{n} \PP\left( \widehat{\bf X}^{(i)}\,\widetilde{\Theta}_i(t)\,{\bf 1}_{\{\theta>0\}} \in x\,A \right) \sim \sum_{i=1}^{n} \PP\left({\bf X}^{(i)}\,e^{-Y_i^*}\,\widetilde{\Theta}_i(t)\,{\bf 1}_{\{\theta>0\}} \in x\,A \right) \,,
\eeam 
uniformly for $n \in \bbn$, and $t \in \Lambda_T$. Taking $n=\infty$ in relation \eqref{eq.KPY.5.32}, we obtain relation \eqref{eq.KPY.5.23}. 
~\halmos

\subsection{Dependence in MRV case} \label{sub.sec.KPT29U.5.3}

Now we focus in two examples from $MRV$ case, the one under weak and the other under strong dependence structure between ${\bf X}^{(i)}$ and $Z_i = e^{-Y_i}$, that are based on the multivariate versions of Breiman theorem, given in 
\cite{fougeres:mercadier:2012}. Our examples provide clear picture of Corollary \ref{cor.KPY.3.1}, about the effect of the dependence between insurance and financial risks on the asymptotic behavior of relation \eqref{eq.KPY.1.7}.

Here we use an equivalent definition of the $MRV(\alpha,\,V,\,\mu)$, with respect to relation \eqref{eq.KPY.2.6}. We say that $F \in MRV(\alpha,\,V,\,\mu)$, if it holds
\beam \label{eq.KPY.5.33}  
\lim x\,\PP\left( \dfrac{{\bf X}}{U_V(x)}\in \bbb \right) = \mu(\bbb) \,,
\eeam 
for a Radon measure $\mu$, non-degenerate to zero, for any Borel set $\bbb \in [0,\,\infty]^d$, with ${\bf 0} \not\in \overline{\bbb}$, $\mu(\partial \bbb)=0$, where
\beao
U_V(x) =\left( \dfrac 1{\bV}\right)^{\leftarrow}(x)\,,
\eeao
which is named normalization function, where $f^{\leftarrow}$ represents a c\'{a}gl\'{a}d inverse of function $f$, that is such that  $U_V \in \mathcal{R}_{1/\alpha}$, see 
\cite[Prop. 2.6]{resnick:2007}.

We also observe that if  $F \in MRV(\alpha,\,V,\,\mu)$, with $\alpha >0$, then by \cite[Th. 1.1(i)]{basrak:davis:mikosch:2002}, we obtain that the $|{\bf X}|$ has regularly varying distribution with the same index $\alpha>0$, hence there exists some non-degenerate to zero, Radon measure $\mu^*$, such that it holds
\beam \label{eq.KPY.5.34}  
\lim x\,\PP\left( \dfrac{|{\bf X}|}{U_V(x)}\in \bbb \right) = \mu^*(\bbb) \,,
\eeam
for any Borel set $\bbb \in [0,\,\infty]$, with $0 \not\in \overline{\bbb}$, $\mu^*(\partial \bbb)=0$.

\subsubsection{Weak dependence scenario} \label{subsubsec.KPT29U.5.3.1}

Further, we present a corollary for the case of asymptotic independence between the ${\bf X}^{(i)}$ and $Z_i = e^{-Y_i}$, inspired by \cite[Th. 2]{fougeres:mercadier:2012}. Hereafter, we consider the inclusion $F \in MRV(\alpha,\,V,\,\mu)$, with $\alpha >0$, except otherwise stated.

\begin{assumption} \label{ass.KPY.5.2}
For any $i \in \bbn$, we suppose that 
\beam \label{eq.KPY.5.35}  
\lim x\,\PP\left( \left[\dfrac{{\bf X}^{(i)}}{U_V(x)},\,Z\right]\in \bbb \right) = (\mu \times L_i)(\bbb) \,,
\eeam 
for any Borel set $\bbb \in [0,\,\infty]^{d+1}$, with ${\bf 0} \not\in \overline{\bbb}$, $(\mu \times L_i)(\partial \bbb)=0$, where for any $i \in \bbn$, the $L_i$ is a proper distribution on $[0,\,\infty)$. Further, we suppose that for all $i \in \bbn$, we have $\E[Z_i^{\alpha}] < \infty$, and there exists some $\delta > 0$, such that it holds
\beam \label{eq.KPY.5.36}  
\lim_{\vep \to 0} \limsup x\,\E \left[\left(\dfrac{|{\bf X}^{(i)}|\,Z_i}{U_V(x)}\right)^{\delta}\,{\bf 1}_{\{|{\bf X}^{(i)}| \leq \vep\,U_V(x)\}}\right] = 0\,.
\eeam 
\end{assumption}

In Assumption \ref{ass.KPY.5.2} the $L_i$ is not necessary the distribution of $Z_i$, for $i \in \bbn$. In some cases of weak dependencies, the $L_i$ can be written as function of the distribution $H_i$ of the the random variable $Z_i$, as we can see in the following example. In this example we use a kind of dependence that is stricter than Assumption \ref{ass.KPY.5.1}. For convenience, we take into consideration only the pair $({\bf X},\,Z)$.

\bexam \label{exam.KPY.5.1}
Let consider that $F \in MRV(\alpha,\,V,\,\mu)$, and there exists some measurable function $h\;:\;[0,\,\infty) \to (0,\,\infty)$ such that it holds
\beam \label{eq.KPY.5.37}  
0 < \inf_{s\in S(Z)} h(s) \leq \sup_{s\in S(Z)} h(s) < \infty\,,
\eeam
and for all Borel sets $\bbb \in[0,\,\infty]^d$, with ${\bf 0} \not\in \overline{\bbb}$, $\mu (\partial \bbb)=0$, we have
\beam \label{eq.KPY.5.38}  
\PP\left( {\bf X} \in x\,\bbb\;|\;Z=z \right) \sim h(z)\,\PP\left( {\bf X} \in x\,\bbb\right) \,,
\eeam
uniformly for $z \in S(Z)$. 

We show that \eqref{eq.KPY.5.35} is true with $L=H^*$, where $H^*(dz)=h(z)\,H(dz)$ represents a  proper distribution on $[0,\,\infty)$, see Remark \ref{rem.KPY.5.1}. Let $\bbb^*:= \bbb \times \bbb'$, where $\bbb \in [0,\,\infty]^d$, as before, and $\bbb' \in [0,\,\infty]$, such that it holds $\PP(Z \in \bbb')>0$ (as in the rest cases \eqref{eq.KPY.5.35} is valid directly as degenerate to zero). Then we obtain
\beao
&&\lim x\,\PP\left( \left[\dfrac{{\bf X}}{U_V(x)},\,Z\right]\in \bbb^* \right) =\lim x\,\PP\left( {\bf X}\in U_V(x)\,\bbb\,,\;Z \in \bbb' \right) \\[2mm] 
&&= \lim \int_{\bbb'} x\,\PP\left( {\bf X}\in U_V(x)\,\bbb\;\big|\;Z =z \right)\,\PP\left( Z \in dz\right) \\[2mm] 
&&= \lim \int_{\bbb'} \dfrac 1{\bV(x)}\,\PP\left( {\bf X}\in x\,\bbb\;\big|\;Z =z \right)\,\PP\left( Z \in dz\right) \\[2mm] 
&&= \lim \dfrac 1{\bV(x)}\, \int_{\bbb'} h(z) \PP\left( {\bf X}\in x\,\bbb\right)\,\PP\left( Z \in dz\right) \\[2mm] 
&&= \lim \dfrac 1{\bV(x)}\, \PP\left( {\bf X}\in x\,\bbb\right)\,\PP\left( Z^* \in  \bbb'\right)=\mu(\bbb)\,H^*(\bbb')=(\mu \times H^*)(\bbb^*)\,,
\eeao
where at the third step we took into consideration the equivalence of \eqref{eq.KPY.5.33} and \eqref{eq.KPY.2.6}, at the pre-last step we used  \eqref{eq.KPY.2.6}, while at the fourth step we employed  \eqref{eq.KPY.5.38}, via its uniformity, since $\bbb' \in S(Z)$. Hence, we find $L=H^*$. We also mention that for the interesting sets $\bbb' \subsetneq S(Z)$, it holds $\mu(\partial \bbb)=0$ if and only if 
\beao
(\mu \times H_i^*)(\partial \bbb^*)=0\,.
\eeao
\eexam 

The following proposition represents a partial reformulation of \cite[Th. 2]{fougeres:mercadier:2012}.

\bpr \label{pr.KPY.5.2}
Let suppose that $({\bf X},\,Z)$ satisfies Assumption \ref{ass.KPY.5.2}. Then the inclusion ${\bf X}\,Z \in MRV(\alpha,\,V,\,\mu_L)$ is valid, that means for any Borel set $\bbb \in [0,\,\infty]^d$, with ${\bf 0} \not\in \overline{\bbb}$ and $\mu_L(\partial \bbb) = 0$, it holds
\beam \label{eq.KPY.5.39}  
\lim \dfrac 1{\bV(x)}\, \PP\left( {\bf X}\,Z \in x\,\bbb \right) =\mu_L\left( \bbb\right) \,,
\eeam
where 
\beam \label{eq.KPY.5.40}  
\mu_L \left( \bbb \right)=(\mu \times L)\left( \{({\bf x},\,y)\;:\; {\bf x}\,y \in \bbb \} \right) =\int_0^{\infty} \mu(y^{-1}\,\bbb)\,L(dy) \,.
\eeam
\epr

In the following corollary we denote for any $i \in \bbn$
\beam \label{eq.KPY.5.41}  
\mu_{L_i}(\bbb) := \int_0^{\infty} \mu(y^{-1}\,\bbb)\,L_i(dy) \,.
\eeam

\bco \label{cor.KPY.5.2}
Under the conditions of Corollary \ref{cor.KPY.3.1}, without $G_i \in MRV$, and under the condition \eqref{eq.KPY.5.35}, it holds
\beam \label{eq.KPY.5.42} 
&&\PP({\bf D}(t) \in x\,A) \\[2mm] \notag
&&\sim \bV(x)\left(\mu(A) \sum_{i=1}^{\infty} \E\left[\Theta_i^{\alpha}(t) {\bf 1}_{\{\theta =0\}} \right]+ \sum_{i=1}^{\infty} \mu_{L_i}(A)\E\left[\widetilde{\Theta}_i^{\beta}(t) {\bf 1}_{\{\theta>0\}} \right] \right) ,
\eeam
uniformly for $t \in \Lambda_T$, for any fixed $T \in \Lambda$. 
\eco

\bre \label{rem.KPY.5.2}
From \eqref{eq.KPY.5.42} we see that the decay rate of probability in \eqref{eq.KPY.1.7} is determined by $V \in \mathcal{R}_{-\alpha}$ (so its magnitude is of order $x^{-\alpha}$), while the rest quantities in the right hand member are positive constants. In the previous case the main role is played by the heaviness of the tails of the insurance risks, while the dependence between the insurance and financial risks appears only as influence on these constants. To see this, we take into account Example \ref{exam.KPY.5.1}, in two cases: in the first case it holds $h_i(z) \neq 1$, for all $z \in S(Z)$, and then we obtain
\beao
\mu_{L_i}(\bbb) = \int_0^{\infty} \mu(y^{-1}\,\bbb)\,H_i^*(dy)= \int_0^{\infty} h_i(y)\,\mu(y^{-1}\,\bbb)\,H_i(dy)\,,
\eeao
while in the second case, when the ${\bf X},\,Z$ are independent each other, then we find 
\beao
h_i(z)=1\,,
\eeao
and hence
\beao
\mu_{L_i}(\bbb) = \int_0^{\infty} \mu(y^{-1}\,\bbb)\,H(dy)\,.
\eeao
\ere

Next, we proceed to the proof of Corollary \ref{cor.KPY.5.2}.

\noindent{\bf Proof of Corollary \ref{cor.KPY.5.2}.}~
From Corollary \ref{cor.KPY.3.1}, we see that it is enough to apply Proposition \ref{pr.KPY.5.2}. For this purpose, we check if \eqref{eq.KPY.5.36} is satisfied.

Let us consider some fixed $\delta > 0$ with $p_1 :=2\,\delta < \alpha$. By 
H\"{o}lder inequality, for all $i \in \bbn$ we obtain
\beao
&&\E\left[\left( \dfrac{|{\bf X}^{(i)}|\,Z_i}{U_V(x)} \right)^{\delta} {\bf 1}_{\{|{\bf X}^{(i)}| \leq \vep\,U_V(x)\}}\right] \\[2mm]
&&=\E\left[\left( \dfrac{|{\bf X}^{(i)}|\,Z_i}{U_V(x)} \right)^{\delta}\;\Big|\; |{\bf X}^{(i)}| \leq \vep\,U_V(x)\right]\,\PP\left( |{\bf X}^{(i)}| \leq \vep\,U_V(x)\right)  \\[2mm]
&&\leq \PP\left( |{\bf X}^{(i)}| \leq \vep\,U_V(x)\right)\,\E^{1/2}\left[\left( \dfrac{|{\bf X}^{(i)}|}{U_V(x)} \right)^{2\,\delta}\;\Big|\; |{\bf X}^{(i)}| \leq \vep\,U_V(x)\right] \\[2mm]
&&\times \E^{1/2}\left[\left( \dfrac{Z_i}{U_V(x)} \right)^{2\,\delta}\;\Big|\; |{\bf X}^{(i)}| \leq \vep\,U_V(x)\right] \leq C\, \PP\left( |{\bf X}^{(i)}| \leq \vep\,U_V(x)\right)\, \left( \dfrac 1{U_V(x)} \right)^{4\,\delta}\,, 
\eeao 
where in the last step, the constant $C> 0$ is derived by the inequalities $M_{p_1}^{(i)} < \infty$ and $\E\left[|{\bf X}^{(i)}|^{2\,\delta}\right] < \infty$, since the $|{\bf X}^{(i)}|$ are regularly varying with index $\alpha > 2\,\delta$. Hence, from the last relation, in combination with \eqref{eq.KPY.5.34} we find that for all $i \in \bbn$ it holds
\beam \label{eq.KPY.5.43}  
&&\lim_{\vep \to 0} \limsup x\,\E\left[\left( \dfrac{|{\bf X}^{(i)}|\,Z_i}{U_V(x)} \right)^{\delta} {\bf 1}_{\{|{\bf X}^{(i)}| \leq \vep\,U_V(x)\}}\right] \\[2mm] \notag
&&\leq C\, \lim_{\vep \to 0} \limsup x\,\PP\left( |{\bf X}^{(i)}| \leq \vep\,U_V(x)\right)\, \left( \dfrac 1{U_V(x)} \right)^{4\,\delta} \\[2mm] \notag
&&= C\, \lim_{\vep \to 0} \limsup \left( \dfrac 1{U_V(x)} \right)^{4\,\delta}
\,\mu^*([0,\,\vep]) \leq C\,C^*\,\lim_{\vep \to 0} \limsup \left( \dfrac 1{U_V(x)} \right)^{4\,\delta}  =0\,,
\eeam
where the constant $C^*>0$ stems from the fact that $\mu^*$ represents a Radon measure. From \eqref{eq.KPY.5.43}  we obtain \eqref{eq.KPY.5.36}. We should also mention that by Proposition \ref{pr.KPY.5.2} is immediately implied that $G_i(x\,A) \asymp G_j(x\,A)$, for any $A \in \mathscr{R}$ and with $i,\,j \in\bbn\,,\; i \neq j$.  
~\halmos

\subsubsection{Strong dependence scenario} \label{subsubsec.KPT29U.5.3.2}

Further, we examine the case of asymptotic independence between the insurance and financial risks ${\bf X}^{(i)}$ and $Z_i$, inspired by \cite[Prop. 7.6]{resnick:2007}. For this we need the assumption of a non-standard $MRV$ structure, as follows.

\begin{assumption} \label{ass.KPY.5.3}
For any $i \in \bbn$, let us suppose that $F \in MRV(\alpha,\,V,\,\mu)$, $H_i \in \mathcal{R}_{-\beta}$, with $\alpha,\,\beta > 0$. We also suppose that there exists some non-degenerate to zero, Radon measure $\nu_i$, such that it holds 
\beam \label{eq.KPY.5.44}  
\lim x\,\PP\left( \left[\dfrac{{\bf X}^{(i)}}{U_V(x)},\,\dfrac{Z_i}{U_{H_i}(x)}\right]\in \bbb \right) = \nu_i(\bbb) \,,
\eeam 
for any Borel set $\bbb \in [0,\,\infty]^{d+1}$, with ${\bf 0} \not\in \bbb$, $\nu_i (\partial \bbb)=0$. Further, we suppose that for some $j =1,\,\ldots,\,d$, we have 
\beam \label{eq.KPY.5.45}  
\limsup x\,\PP\left[\left(\dfrac{ X_j^{(i)}}{U_V(x)},\,\dfrac{Z_i}{U_{H_i}(x)}\right)\in (0,\,\infty]^2\right] = \nu_i((0,\,\infty]^2) >0\,.
\eeam 
\end{assumption}

\bre \label{rem.KPY.5.3}
Relation \eqref{eq.KPY.5.45} implies the asymptotic dependence of $Z_i =e^{-Y_i}$ with at least one of the components of ${\bf X}^{(i)}$.
\ere

The following proposition is a reformulation of \cite[Prop. 7.6]{resnick:2007}, see also \cite[Th. 4]{fougeres:mercadier:2012}, for an extension of the result to Hadamard products. In the following proposition, for sake of convenience, we use the pair $({\bf X},\,Z)$. We also denote $\overline{Q_i}(x):=[1/(U_V\,U_{H_i})^{\leftarrow}(x)] \wedge 1$, with $i \in \bbn$, for the distribution $Q_i \in \mathcal{R}_{-\alpha\,\beta/(\alpha+ \beta)}$.

\bpr \label{pr.KPY.5.3}
Under Assumption \ref{ass.KPY.5.3} it holds ${\bf X}\,Z \in MRV\left(\dfrac{\alpha\,\beta}{\alpha + \beta},\,Q,\,\nu \right)$, with
\beam \label{eq.KPY.5.46}  
\PP\left({\bf X}\,Z \in x\,\bbb\right) \sim \dfrac{\nu(\bbb^*)}{(U_V\,U_H)^{\leftarrow}(x)} \sim \nu(\bbb^*)\,\overline{Q}(x)\,,
\eeam
for any Borel set $\bbb \in [0,\,\infty]^{d}$, with ${\bf 0} \not\in \bbb$, $\nu (\partial \bbb)=0$, where $\bbb^*:=\{({\bf x},\,z) \in [0,\,\infty)^d \times [0,\,\infty)\;:\; {\bf x}\,z \in \bbb \}$.
\epr

In the following corollary, we demand that $H_i \in \mathcal{R}_{-\beta}$, with the same index of regular variation. The reason is that if $\beta_1 \neq \beta_2$, and $H_i \in \mathcal{R}_{-\beta_i}$, for $i =1,\,2$, then under the conditions of Proposition \ref{pr.KPY.5.3} we find that $G_i \in  MRV\left(\dfrac{\alpha\,\beta_i}{\alpha + \beta_i},\,Q_i,\,\nu \right)$, for $i =1,\,2$, and hence the condition  $G_1(x\,A) \asymp G_2(x\,A)$ fails. Additionally we need that $\beta > \alpha$, see Remark \ref{rem.KPY.3.1}. The proof of the following corollary follows directly from Corollary \ref{cor.KPY.3.1}, Proposition \ref{pr.KPY.5.3} and the fact that the infinite sums at the right member of relation \eqref{eq.KPY.3.10} are finite.

\bco \label{cor.KPY.5.3}
Let suppose that conditions of Corollary \ref{cor.KPY.3.1} are satisfied, without the $G_i \in MRV$. Under Assumption \ref{ass.KPY.5.3} with $\alpha < \beta$ it holds
\beam \label{eq.KPY.5.47}  
\PP\left({\bf D}(t) \in x\,A \right) \sim\sum_{i=1}^{\infty}  \overline{Q}_i(x)\,\nu_i(A^*) \E\left[\widetilde{\Theta}_i^{\beta}(t) {\bf 1}_{\{\theta>0\}} \right]\,,
\eeam
uniformly for $t \in \Lambda_T$, for any fixed $t \in \Lambda_T$, where
\beao
A^*:=\{({\bf x},\,z) \in [0,\,\infty)^d \times [0,\,\infty)\;:\; {\bf x}\,z \in A \}\,.
\eeao
\eco

\bre \label{rem.KPY.5.4}
From \eqref{eq.KPY.5.47} we see that the decay rate of the probability from \eqref{eq.KPY.1.7} is drived by a distribution $Q_i \in  \mathcal{R}_{-\alpha\,\beta/(\alpha+ \beta)}$, while the infinite sum in the right member of \eqref{eq.KPY.5.47} represents only a constant. Since $\alpha\,\beta/(\alpha+ \beta) < \alpha$, we obtain that the decay rate here is slower in comparison with the case of asymptotic independence from Corollary \ref{cor.KPY.5.2}, recall that there it was found with order of magnitude $x^{-\alpha}$. Hence, it indicates the practical importance of the role of dependence between the insurance and financial risks.
\ere

\subsection{Examples on Nested Archimedean copulas} \label{sec.KPT29U.6}

In this last section, we employ two-level nested Archimedean copulas to model the dependence between the non-negative claim vector \({\bf X}=(X_1,\ldots,X_d)\) and the non-negative factor \(Z=e^{-Y}\), in the case of $MRV$, in order to stress the generality of our main results. Our purpose is 
to derive the corresponding nonstandard multivariate regular-variation measure of \(({\bf X},Z)\). We introduce a nested Gumbel copula to model upper-tail asymptotic dependence and a nested Frank 
copula to model upper-tail asymptotic independence.

The $d$-dimensional exchangeable Archimedean construction generated by $\psi$ is
\beam	\label{eq:exchangeable-archimedean}
	C_{\psi}(u_1,\ldots,u_d) 	=	\psi\left(\sum_{i=1}^{d}\psi^{-1}(u_i)\right),
	\qquad (u_1,\ldots,u_d)\in[0,1]^d.
\eeam
where $\psi \colon [0,\infty]\to [0,1]$ is any continuous, decreasing, convex function satisfying $\psi(0)=1$ and $\psi(\infty):=\lim_{t\to\infty}\psi(t)=0$. The function $\psi$ is known as the copula generator and $\psi^{-1}$ is its 
c\'{a}gl\'{a}d inverse, defined in general by $\psi^{-1}(u)=\inf\{t : \psi(t)=u\}$. See \cite{nelsen:2006} for a textbook treatment of copulas.

A sufficient and necessary condition for \eqref{eq:exchangeable-archimedean} to define a copula in every dimension $d\geq2$ is that $\psi$ be completely monotone on
$[0,\infty)$; that is, $\psi$ is continuous at zero, infinitely differentiable on $(0,\infty)$, and satisfies
\beao
(-1)^k \dfrac{d^k}{dt^k}\psi(t)\ge 0,\quad k\in\mathbb{N},\ t>0.
\eeao
This result is due to \cite{kimberling:1974}; see also \cite{mcneil:2008} and  \cite{mcneil:neslehova:2009}.

Let $U_j=F_j(X_j),\ j=1,\ldots,d,\ V=H(Z)$.
Continuity of the marginal distributions implies that all the variables \(U_1,\ldots,U_d,V\) are standard uniform.  We now consider the two-level nested Archimedean copula
\beam	\label{eq:general-nested}
C({\bf u},v)
	=\psi_0\left[	\psi_0^{-1}\{C_1({\bf u})\}+\psi_0^{-1}(v)\right]=\psi_0\left[\psi_0^{-1}\left\lbrace \psi_1\left(\sum_{j=1}^d\psi_1^{-1}(u_j)\right)\right\rbrace +\psi_0^{-1}(v)\right].
\eeam
Here \(v\) corresponds to the outer variable $V=H(Z)$, while $C_1({\bf u})$ is the inner block. It should be emphasized that the nesting construction does not imply that the random variable $C_1(U_1,\ldots,U_d)$ is uniformly distributed. Rather, the validity of \eqref{eq:general-nested} as a $(d+1)$-dimensional copula follows 
from \cite[Th. 4.4]{mcneil:2008}, provided that $\psi_0^{-1}\circ\psi_1$ has a completely monotone derivative. Thus, its validity is established by the nested Archimedean copula theorem, rather than by applying Sklar's theorem directly to $C_1(U_1,\ldots,U_d)$.

\subsubsection{Asymptotic dependence: the nested Gumbel model.} \label{subsec.KPT29U.6.1}

We first consider a nested Gumbel--Hougaard copula. For
\(\theta\geq1\), the completely monotone Gumbel generator is
\beao
\psi_{\theta}^{G}(t)= \exp\left(-t^{1/\theta}\right),
\qquad t\geq0,
\eeao
with inverse
\beao
\left(\psi_{\theta}^{G}\right)^{-1}(u)= (-\log u)^\theta, \qquad u\in(0,1].
\eeao
The Gumbel copula is upper-tail asymptotically dependent whenever
\(\theta>1\); see, for example \cite[Th. 4.1]{charpentier:segers:2009}.

Let \(\theta_0\) be the parameter of the outer copula connecting
\({\bf X}\) and \(Z\), and let \(\theta_1\) be the parameter
of the inner copula describing the dependence among
\(X_1,\ldots,X_d\). Define
\beao
C_1^G({\bf u})=\exp\left\{-\left(\sum_{j=1}^d(-\log u_j)^{\theta_1}\right)^{1/\theta_1}\right\}.
\eeao
The corresponding two-level nested copula is
\begin{equation}\label{eq:nested-gumbel}
		C_G({\bf u},v)
		=C_{\theta_0}^G\left(C_1^G({\bf u}),v\right)
		=\exp\left\{	-\left[	\left(	\sum_{j=1}^d(-\log u_j)^{\theta_1}	\right)^{\theta_0/\theta_1}+(-\log v)^{\theta_0}\right]^{1/\theta_0}\right\},
\end{equation}
where \(1<\theta_0\leq\theta_1<\infty\). Under this parameter restriction, the sufficient nesting condition of \cite{joe:1997} and \cite{mcneil:2008} is satisfied. Hence, \(C_G\) defines a valid \((d+1)\)-dimensional nested Archimedean copula.

The following result gives the specific form of the limit measure of the nonstandard multivariate regular variation of \(({\bf X},Z)\) under nested Gumbel copula. This result can be viewed as a generalization 
of \cite[Lem. 5.2]{tang:yuan:2013}.

\begin{example} \label{exam.KPY.a}
Suppose that the non-negative random vector \({\bf X}=(X_1,\ldots,X_d)\) has continuous marginal distribution functions \(F_1,\ldots,F_d\) satisfying satisfy the tail equivalence condition
	$\overline{F}_j(x)\sim c_j\overline{F^*}(x),\ j=1,\ldots,d,$	where \(c_1,\ldots,c_d\) are positive constants and	\( F^* \in\mathcal{R}_{-\alpha}\) for some \(\alpha \in (0,\,\infty)\).
	Let \(Z\) have distribution function \(H\) satisfying	\(H \in\mathcal{R}_{-\beta}\) for some \(\beta \in (0,\,\infty)\).
	Suppose further that \(({\bf X},Z)\) possesses the nested Gumbel copula \eqref{eq:nested-gumbel}, with
	\(1<\theta_0\leq\theta_1<\infty\). Then, 
\beao
\lim x\,\PP\left[
	\left(\dfrac{\bf X}{U_{F^*}(x)},\dfrac{Z}{U_H(x)}\right)\in\mathbb{B}\right] = \nu_G(\mathbb{B}) \,,
\eeao
for any Borel set $\mathbb{B} \subset [0,\,\infty]^{d+1}$, with ${\bf 0} \not\in \overline{\mathbb{B}}$, $\nu_G (\partial \mathbb{B})=0$, where
the limit measure \(\nu_G\) is determined by
	\begin{equation} \label{eq:gumbel-exponent-measure}
		\nu_G\left(\left([0,{\bf t}]\times[0,z]\right)^c\right)=
		\left[\left(\sum_{j=1}^d c_j^{\theta_1}\,t_j^{-\alpha\theta_1}\right)^{\theta_0/\theta_1}+z^{-\beta\theta_0}\right]^{1/\theta_0},
	\end{equation}
	for every \({\bf t}\in(0,\infty)^d\) and \(z>0\), and $U_{F^*},\,U_H$ represent the normalization functions, corresponding to distributions $F^*$ and $H$.
\end{example}

\begin{proof}
Let ${\bf U}=(U_1,\ldots,U_d,V)$ be a \((d+1)\)-dimensional uniform random vector distributed according to the nested Gumbel copula \(C_G\), defined by relation \eqref{eq:nested-gumbel}. For every ${\bf a}=(a_1,\ldots,a_d,a_z)\in(0,\infty)^{d+1}$, by (\ref{eq:nested-gumbel}) and Taylor expansion we have
\begin{align*}
&\lim x\,\PP\left(\bigcup_{j=1}^d\left\{U_j>1-\dfrac{a_j}{x}\right\}\cup\left\{V>1-\dfrac{a_z}{x}\right\}\right) 		\\[2mm]
&=\lim x\,\left[1-C_G\left(1-\dfrac{a_1}{x},\ldots,1-\dfrac{a_d}{x},1-\dfrac{a_z}{x}\right)\right]=\left[\left(\sum_{j=1}^d a_j^{\theta_1}\right)^{\theta_0/\theta_1}+a_z^{\theta_0}\right]^{1/\theta_0}.
\end{align*}
Then we have the upper bound
\begin{align*}
	&x\,\PP\left[\left(\dfrac{\bf X}{U_{F^*}(x)},\dfrac{Z}{U_H(x)}\right)\in\left([{\bf 0},{\bf t}]\times[0,z]\right)^c\right]\\[2mm]
	&=x\,\left[	1-C_G\left(	F_1\left(U_{F^*}(x)\,t_1\right),\ldots,F_d\bigl(U_{F^*}(x)\,t_d\bigr),H\bigl(U_H(x)\,z\bigr)\right)\right]\\[2mm]
	&=x\,\PP\left(\bigcup_{j=1}^d\left\{U_j>1-\overline{F}_j\bigl(U_{F^*}(x)\,t_j\bigr)\right\}\cup\left\{V>1-\overline{H}\bigl(U_H(x)\,z\bigr)\right\}\right)\\[2mm]
	&\lesssim x\,\PP\left(\bigcup_{j=1}^d\left\{U_j>1-\dfrac{(1+\varepsilon)c_j\,t_j^{-\alpha}}{x}\right\}\cup\left\{V>1-\dfrac{(1+\varepsilon)z^{-\beta}}{x}\right\}\right)\\[2mm]
	&\rightarrow (1+\varepsilon)\,\left[\left(\sum_{j=1}^d c_j^{\theta_1}\,t_j^{-\alpha\,\theta_1}\right)^{\theta_0/\theta_1}+z^{-\beta\theta_0}\right]^{1/\theta_0}.
\end{align*}
Because of the arbitrary choice of $\vep >0$, we obtain the desired upper bound for \eqref{eq:gumbel-exponent-measure}.

Establishing the corresponding lower bound by a similar way, completes the proof.
\end{proof}

We should mention that under the conditions of the previous example, from \cite[Lem. 5.2]{tang:yuan:2013}, we find that ${\bf X} \stackrel{d}{\sim} F \in MRV(\alpha,\,F^*,\,\mu)$, 
where the Radom measure $\mu$ is of the form
\beao
\mu\left( [{\bf 0},\,{\bf t}]^{c}\right) =\left( \sum_{j=1}^{d} c_j^{\theta_1}\,t_j^{-\alpha\,\theta_1} \right)^{\theta_0 / \theta_1}\,,
\eeao
with ${\bf t} \in ({\bf 0},\,{\bf \infty})^{d}$. This can be also checked following the line of the proof of Example \ref{exam.KPY.a}, only for ${\bf X}$. Hence, Example \ref{exam.KPY.a}, 
 when $\alpha < \beta$, satisfies all the conditions of Proposition \ref{pr.KPY.5.3}. 

Further, when \(d=1\), formula \eqref{eq:gumbel-exponent-measure} reduces to
\beao
\nu_G([0,t]\times[0,z])=\left(c_1^{\theta_0}\,t^{-\alpha\theta_0}+z^{-\beta\theta_0}\right)^{1/\theta_0}, 
\eeao 
which is consistent with the measure result for the case when $(X,Z)$ follows $BRV$ in \cite[Cor. 3.2]{yang:fan:yuen:2023}.

\subsubsection{Asymptotic independence: the nested Frank model} \label{subsec.KPT29U.6.2}

We next construct an upper-tail asymptotically independent copula model and we provide a sufficient condition, such that it is satisfied Example \ref{exam.KPY.5.1}, and further by Corollary \ref{pr.KPY.5.2}.
For $\theta>0$, the Frank generator is
\beao
\psi_\theta^F(t)=-\dfrac1\theta\log\left\{1-(1-e^{-\theta})e^{-t}\right\},\qquad t\geq0,
\eeao
with inverse
\beao
(\psi_\theta^F)^{-1}(u)=-\log\left(\dfrac{1-e^{-\theta u}}{1-e^{-\theta}}\right),\qquad 0<u\leq1.
\eeao
Let \(\theta_0\) be the parameter of the outer copula connecting \({\bf X}\) and \(Z\), and let \(\theta_1\) be the parameter of the inner copula describing the dependence among \(X_1,\ldots,X_d\). 
Let the inner copula be
\beao
C_1^F({\bf u})=-\dfrac1{\theta_1}\log\left[1-\dfrac{\prod_{j=1}^d(1-e^{-\theta_1u_j})}{(1-e^{-\theta_1})^{d-1}}\right].
\eeao
The corresponding two-level nested copula is
\beam \label{eq.KPY.6.A}
C_F({\bf u},v)=C_{\theta_0}^{F}\bigl(C_1^F({\bf u}),v\bigr)
	=-\dfrac1{\theta_0}\log\left[1-\dfrac{\left( 1-e^{-\theta_0 C_1^F({\bf u})}\right) (1-e^{-\theta_0v})}{1-e^{-\theta_0}}\right],
\eeam
where $0<\theta_0\leq\theta_1<\infty$ is the sufficient nesting condition for a valid $(d+1)$-dimensional nested Archimedean copula.

Retain the strong tail equivalence condition, (namely $\bF_j(x) \sim  c_j\,\bF^*(x)$, for $j=1,
,\ldots,\,d$ with $(c_1,\,\ldots,\,c_d)\in(0,\,\infty)^{d}$ and $F^* \in \mathcal{R}_{-\alpha}$, with $\alpha \in (0,\,\infty)$.
). Since the inner Frank copula is upper-tail asymptotically independent, the 
$MRV$ limit measure of  ${\bf X}$ is concentrated on the coordinate axes. In particular,
\beao 
	\mu([{\bf 0},{\bf t}]^c)=\sum_{j=1}^dc_jt_j^{-\alpha},\qquad {\bf t}\in(0,\infty)^d.
\eeao

\bexam \label{exam.KPY.6.f}
Let us consider that the $({\bf X},\,Z)$ satisfies the nested copula condition \eqref{eq.KPY.6.A}. 
Further we suppose that the marginal distributions of ${\bf X}$, are such that it holds 
$\bF_j(x) \sim c_j\,\bF^*(x)$ for some positive constants $c_j$, for $j=1,\,\ldots,\,d$, and 
$F^* \in \mathcal{R}_{-\alpha}$, with $\alpha \in (0,\,\infty)$. Then \eqref{eq.KPY.5.38} holds uniformly 
for $z \in S(Z)$.
\eexam

\pr~
We next proceed to verify the assumption under the nested Frank model for the given set $B_{\bf t}:=[{\bf 0},\,{\bf t}]^c,\ {\bf t}:=(t_1,\ldots,t_d)\in(0,\infty)^d$.
Direct differentiation of the nested Frank copula gives
\beao
\dfrac{\partial C_F({\bf u},v)}{\partial v}=\dfrac{e^{-\theta_0v}\left(1-e^{-\theta_0C_1^F({\bf u})}\right)}{(1-e^{-\theta_0})-\left(1-e^{-\theta_0C_1^F({\bf u})}\right)\left(1-e^{-\theta_0v}\right)}.
\eeao
Then
\begin{align*}
	\PP({\bf X} \in xB_{\bf t}\;\big|\;  Z=z)	&=
	1-\dfrac{\partial C_F({\bf u},v)}{\partial v}
	\big|_{{\bf u}={\bf u}_x,\,v=H(z)}\\
	&=
	\dfrac{e^{-\theta_0C_1^F({\bf u}_x)} -e^{-\theta_0} }{(1-e^{-\theta_0})-	\left( 1-e^{-\theta_0 C_1^F({\bf u}_x} )\right) 
		\left( 1-e^{-\theta_0 H(z)}\right) }.
\end{align*}
where ${\bf u}_x=\bigl(F_1(xt_1),\ldots,F_d(xt_d)\bigr),\ v=H(z)$. It follows that
\beao
\dfrac{\PP({\bf X}\in xB_{\bf t}\;\big|\; Z=z)}{\PP({\bf X}\in xB_{\bf t})}
=\dfrac{e^{-\theta_0C_1^F({\bf u}_x)}-e^{-\theta_0}}{\left( 1-C_1^F({\bf u}_x)\right) \left[(1-e^{-\theta_0})-\left( 1-e^{-\theta_0C_1^F({\bf u}_x)}\right) \left( 1-e^{-\theta_0H(z)}\right) \right]}.
\eeao
	Hence,
\begin{align*}
	h(z)=	\lim\dfrac{\PP\bigl({\bf X}\in xB_{\bf t}\;\big|\; Z=z\bigr)}{\PP\bigl({\bf X}\in xB_{\bf t}\bigr)}
	=\dfrac{\theta_0 e^{\theta_0 H(z)}}{e^{\theta_0}-1}.
\end{align*}
Since $(1-e^{-\theta_0})e^{-\theta_0H(z)}\geq(1-e^{-\theta_0})e^{-\theta_0}>0$ and $\inf_{z\in S(Z)}h(z)
\geq\dfrac{\theta_0}{e^{\theta_0}-1}>0$  the preceding convergence is uniform for \(z\in S(Z)\), namely it holds
\beao
\PP({\bf X}\in x\,B_{\bf t}\;\big|\; Z=z) \sim h(z)\PP({\bf X}\in x\,B_{\bf t})\,,
\eeao
uniformly for $z \in S(Z)$. Applying in the last relation the Dynkin's $\pi - \lambda$ theorem, 
see \cite[Cor. 2.2.1]{resnick:1999}, we obtain that relation \eqref{eq.KPY.5.38} is true for any $\bbb \in [0,\,\infty]^d$, with ${\bf 0} \not\in \overline{\bbb}$, with the desired 
uniformity with respect to  \(z\in S(Z)\).
~\halmos

\vspace{1em}
\noindent\textbf{Funding Statement}: This work is supported by the Basic Scientific Research Project of the Department of Education of Liaoning Province (No. LJ212510173010) and Doctoral Research Start-up Program of the Department of Science and Technology of Liaoning Province (No. 2026-BS-0753).

\end{document}